\documentclass[11pt]{amsart}
\usepackage{amsmath,amssymb,mathrsfs,color}
\usepackage[title]{appendix}
\usepackage{amssymb}
\usepackage{hyperref}
\allowdisplaybreaks

\let\cal=\mathcal

\def\Z{{\mathbb Z}}
\def\R{{\mathbb R}}
\def\P{{\mathbb P}}
\def\E{{\mathbb E}}

\def\eps{\epsilon}

\def\Var{\operatorname{Var}}
\newtheorem{thm}{Theorem}[section]
\newtheorem{cor}[thm]{Corollary}
\newtheorem{lem}[thm]{Lemma}
\newtheorem{prop}[thm]{Proposition}

\theoremstyle{definition}
\newtheorem{de}[thm]{Definition}
\theoremstyle{remark}
\newtheorem{rem}[thm]{Remark}

\numberwithin{equation}{section}

\newcommand{\rmd}{{\rm d}}

\newcommand{\Lip}{{\operatorname{Lip}}}
\newcommand{\fcb}{{\operatorname{fcb}}}
\newcommand{\Leb}{\operatorname{Leb}}

\allowdisplaybreaks

\begin{document}

\title[Rates in the WIP for slowly mixing sequential dynamical systems]{Rates in the weak invariance principle for slowly mixing sequential dynamical systems}

\author{Nicholas Fleming-V\'{a}zquez}
\address{N. Fleming-V\'{a}zquez: School of Mathematics and Statistics, UNSW Sydney,
Sydney NSW 2052, Australia}
\email{nicholas.fleming@unsw.edu.au}

\author{Zhe Wang}
\address{Z. Wang: School of Mathematical Sciences, Dalian University of Technology, Dalian
116024, P. R. China}
 \email{zwangmath@hotmail.com}

\date{\today}



\begin{abstract}
We obtain rates of convergence in the weak invariance principle for certain classes of slowly mixing sequential dynamical systems. 
In particular, for systems satisfying a Functional Correlation Bound with rate $n^{-\gamma} \ (\gamma>1)$, we establish an explicit bound in the Wasserstein-$r$ metric for $1\le r<\gamma+1$ whenever the variance grows at a sufficiently fast polynomial rate.
The proof is based on Bernstein's ``big block--small block" method.
Examples covered by our result include sequential heterochaos baker maps and  sequential LSV maps. 
In particular, for random LSV maps over an ergodic base, we obtain rates in the quenched weak invariance principle.
\end{abstract}

\maketitle

\section{Introduction}
\setcounter{equation}{0}

Let $T:M\to M$ be a probability-measure-preserving chaotic dynamical system and let
$v:M\to\R$ be a sufficiently regular observable with mean zero. 
The \emph{central limit theorem (CLT)} states that the normalized Birkhoff sums converge in distribution to a centered Gaussian random variable.  
The \emph{weak invariance principle (WIP)}, also known as the functional CLT,
strengthens the CLT by stating that the corresponding partial-sum processes converge in distribution to Brownian motion.
Over the past few decades, the WIP has been established for a wide range of systems, including uniformly and nonuniformly expanding/hyperbolic maps and flows~\cite{hofbauer_keller,denker_philipp,MN05}.

In many applications, the evolution of a system is influenced by parameters or external noise,
so that the map defining the dynamics
may vary from one time step to the next. 
This motivates the study of compositions of a sequence of maps, generalizing the classical setting in which a single transformation is iterated. 
More precisely, given a probability space $(M,m)$ and a sequence of maps $T_n:M\to M$, the compositions
\[
\mathcal T^n=T_n\circ\cdots\circ T_1
\]
define a sequential dynamical system. 
One is then interested in studying limit theorems for partial sums of the form $S_n = \sum_{i=0}^{n-1}v_i \circ \cal T^i$, where $v_i:M\to \R$ are sufficiently regular observables.
Limit theorems for sequential dynamical systems have been developed in~\cite{Conze07,HNTV17,NTV18,Su22,DH25,DL25} and the references therein. 
An important special case is that of random dynamical systems,
where the sequence of maps is generated by a random environment.
This setting has been a major focus of research on limit theorems;
see for example~\cite{K98,ANV15,DFGV18,DH21,Su23}.
Note that the sequence of partial sums $S_n$ is generally nonstationary and its variance need not be asymptotically linear. 
These features make quantitative limit theorems substantially more delicate than in
the autonomous setting. 

Recently, the problem of quantifying the convergence rate in the WIP for dynamical systems has attracted much interest~\cite{antoniou_melbourne,hafouta_alpha_mixing,LW24,F25,paviato}. 
To our knowledge, only two existing results~\cite{LW24,liu_saussol_vaienti_wang} give rates in the Wasserstein metric for sequential dynamical systems.
In \cite{LW24}, Liu and Wang obtained a rate of convergence in the WIP for sequential dynamical systems under a spectral gap type assumption for nonstationary compositions of the transfer operators. 
However, many interesting systems, such as sequential intermittent maps, do not satisfy this assumption. 
In~\cite{liu_saussol_vaienti_wang}, Liu, Saussol, Vaienti and Wang obtained rates in the quenched WIP for random dynamical systems modeled by random Young towers, including random intermittent maps. 
However, they use a ``secondary-martingale coboundary decomposition" which is obtained by exploiting the tower structure; for general sequential intermittent maps it is not clear how to obtain such a decomposition.

The Functional Correlation Bound, introduced by Lepp\"anen in~\cite{L17}, provides a more flexible framework for obtaining quantitative limit theorems.
This condition is satisfied by sequential intermittent maps~\cite{L17}; in the autonomous setting it holds for broad classes of slowly mixing nonuniformly hyperbolic systems~\cite{F22} and leads to Wasserstein rates in the WIP~\cite{F25}. 
It has also been used by Lepp\"anen and coauthors to prove quantitative limit theorems for sequential dynamical systems via Stein's method~\cite{L17,leppanen_stenlund-sunklodas,hafouta_leppanen}, including convergence rates in the WIP in a smooth metric~\cite{leppanen_nakajima_nakano}.

In this paper, we obtain convergence rates in the WIP for sequential dynamical systems that satisfy the Functional Correlation Bound with rate $n^{-\gamma}$, $\gamma>1$. 
Here, we employ the Wasserstein distance to measure the rate of weak convergence.  
For $p\ge 1$, we denote by $\cal W_p(P,Q)$ the Wasserstein-$p$ metric between distributions $P$ and $Q$ on a metric space $(\cal X,d)$ (see \cite[Section~3.3]{bogachev_weak_conv}):
\[
\cal W_p(P,Q)= \inf \big\{ \E [{d(X,Y)}^p]^{1/p}\colon \operatorname{law} (X)=P,\ \operatorname{law} (Y)=Q \big\}.
\]

Our main result, Theorem~\ref{thnon}, gives an explicit convergence rate in the WIP with respect to the Wasserstein-$p$ metric. 
In contrast to many results on limit theorems for sequential dynamical systems~\cite{HNTV17,Su22,leppanen_nakajima_nakano}, we cover examples where the variance $\Sigma_n^2$ grows at an arbitrarily slow polynomial rate. 
More precisely, if the Functional Correlation Bound holds with rate $n^{-\gamma}$, $\gamma>1$, then we obtain rates in the WIP whenever $\lim_{n\to\infty}\Sigma_n/n^\kappa=\infty$ for some $\kappa>1/(\gamma+1)$.
We also provide an example to show that the assumption $\kappa>1/(\gamma+1)$ is essentially sharp in the Appendix.

The proof of our main result is based on Bernstein's classical ``big block--small block" method, as in~\cite{F25}; 
however, various adaptations are necessary to handle our nonstationary setting. 
In~\cite{F25}, the sums over big blocks are identically distributed; here we define the blocks so that the sums over the big blocks have approximately equal variance. 
We also prove a maximal inequality in $L^{\gamma+1}$ for partial sums. 
In particular, when $\Sigma_n^2 =O(n^\theta)$ with $\theta<1$, the bound we obtain is sharper than the usual $O(n^{1/2})$ bound.

Our result covers a large class of sequential dynamical systems, including sequential heterochaos baker maps and sequential LSV maps. 
In particular, for random LSV maps over an ergodic base, it yields convergence rates in the quenched WIP. 
The result also applies in the autonomous setting to mixing nonuniformly hyperbolic maps modeled by Young towers with $O(n^{-(\gamma+1)})$ tails for $\gamma>1$, including Anosov diffeomorphisms, the Sinai billiard map and Bunimovich flowers. 
We discuss these examples in more detail in Section~\ref{sec:examples}.

The remainder of the paper is organized as follows. In Section~\ref{sec:main_result}, we introduce
the Functional Correlation Bound and state the main theorem. 
In Section~\ref{sec:examples}, we provide several examples that our result applies to.
Section~\ref{sec:moment_estimates} establishes the maximal moment estimates needed in the proof of the main result. 
In Section~\ref{sec:main_thm_pf}, we prove the main result. In the \hyperref[sec:variance_growth_sharp]{Appendix}, we give an example to show that the polynomial variance growth assumption in our main result is essentially sharp.

\emph{Notation:} Throughout this paper, $\|\cdot\|_p$ denotes the $L^p$-norm. 
For random variables $X,Y$, we write $X=_d Y$ if $X$ and $Y$ have the same distribution. 
For an integrable random variable $X$ on $(M,m)$, we write $\mathbb E_m[X]$ for the expectation of $X$; for random variables on auxiliary probability spaces, we write $\mathbb E$ for expectation with respect to the underlying probability measure, which is understood from the context.
For random elements $X$, $Y$ of $D[0,1]$ we write
$\mathcal W_p(X,Y)$ for the Wasserstein-$p$ distance between their laws. 
Given nonnegative sequences $a_n$ and $b_n$, we write $a_n\lesssim b_n$ or $a_n = O(b_n)$ if there exists a constant $C>0$ such that $a_n\le Cb_n$ for all $n\ge 1$. In particular, we write $a_n\asymp b_n$ if $a_n\lesssim b_n$ and $b_n \lesssim a_n$.

\section{Setting and main result}\label{sec:main_result}
In this section, we first introduce the Functional Correlation Bound and some of its properties, and then we state our main result.

\subsection{The Functional Correlation Bound}
 
Let $(M,d)$ be a metric space. For $k\ge 1$, we call a bounded function $G:M^k\to \R$ \emph{separately Lipschitz} if $|G|_\infty = \sup |G|<\infty$ and $\Lip(G;i)<\infty$ for all $1\le i\le k$, where 
\begin{equation*}
    \Lip(G;i)=\sup_{x_1,\dots,x_k,\ x'_i\ne x_i\in M}\frac{|G(x_1,\dots,x_k)-G(x_1,\dots,x_{i-1},x'_i,x_{i+1},\dots,x_k)|}{d(x_i,x'_i)}.
\end{equation*}
Let $\{T_n\}_{n\ge 1}$ be a sequence of measurable transformations on $M$ and $m$ be a probability measure on $M$. 
  For $1\le j \le n$, we set
\[
\cal{T}_j^{n}:=T_n\circ T_{n-1}\circ\cdots\circ T_j, \ \ \cal{T}^n:=\cal{T}_1^n, \ \ \cal T^0:=Id.
\]
\begin{de}\label{de:fcb}
	We say that $\{\cal T^n\}_{n\ge 0}$ satisfies the {\em Functional Correlation Bound} with rate $n^{-\gamma}$ (and constant $C_{\fcb}>0$) if for all integers $1\le p<q$ and $0\le n_1\le \dots \le n_q$ and any separately Lipschitz function $G:M^q\to \R$, we have 
	\begin{align}
		& \bigg|\int_M G(\cal T^{n_1}x,\dots,\cal T^{n_q}x)\rmd m(x)\notag \\
		&\qquad \qquad -\iint_{M^2} G(\cal T^{n_1}x_0,\dots, \cal T^{n_p}x_0,\cal T^{n_{p+1}}x_1,\dots, \cal T^{n_q}x_1) \rmd m(x_0)\rmd m(x_1)\bigg|\notag
        \\
		& \qquad \le C_{\fcb}\, (n_{p+1} - n_p)^{-\gamma}\bigg(|G|_\infty + \sum_{i=1}^{q}\Lip(G;i)\bigg).\label{eq:fcb}
	\end{align}
\end{de}
\begin{rem}\label{rem:fcb}
    (i) A condition similar to Definition~\ref{de:fcb} was introduced by Lepp\"anen in~\cite{L17}, where it was verified for sequential intermittent maps. 
    In Section~\ref{sec:examples}, we discuss this example in more detail and give further examples that satisfy this condition.
    
	(ii) Let $0= k_0<k_1<\dots$ be integers. Note that if $\{\cal T^n\}_{n\ge 0}$ satisfies the Functional Correlation Bound with rate $n^{-\gamma}$ and constant $C_{\fcb}$, then $\{\cal T^{k_n}\}_{n\ge 0} $ also satisfies the Functional Correlation Bound with rate $n^{-\gamma}$ and constant $C_{\fcb}$.
    
	(iii) In~\cite{F25}, a Functional Correlation Bound is considered for separately H\"older functions. There is no loss of generality in only considering separately Lipschitz functions: if $(M,d)$ is a metric space, then $(M,d^\eta)$ is also a metric space for any $\eta\in (0,1)$ and (separately) $\eta$-H\"older functions with respect to $d$ are (separately) Lipschitz with respect to $d^\eta$.
\end{rem}

The Functional Correlation Bound can be seen as a weak dependence condition. 
Suppose that $\{\cal T^n\}_{n\ge 0}$ satisfies the Functional Correlation Bound with rate $n^{-\gamma}$ and constant $C_{\fcb}>0$.
Let $k\ge 1$ and consider $k$ disjoint blocks of integers $\{a_i,a_i+1,\ldots,b_i\}$, $1\le i\le k$ with $a_i\le b_i<a_{i+1}$. 
Consider random variables $X_i$ of the form
\[
X_i(x)=\Phi_i(\cal T^{a_{i}}x,\ldots,\cal T^{b_i}x),
\]
where $\Phi_i:M^{b_i - a_i+1}\to \R$ is a separately Lipschitz function. 
Let $\widehat X_1,\dots,\widehat X_{k}$ be independent random variables such that $\widehat X_i =_d X_i$.

\begin{prop}\label{lip}
Let $R=\max_{1\le i\le k}|\Phi_i|_\infty$. Then for any Lipschitz function $\varphi:[-R,R]^k\to \R$, we have
\begin{align*}
&\big|\E_m[\varphi(X_1,\ldots,X_k)]-\E[\varphi(\widehat X_1,\ldots,\widehat X_k)]\big|\\
\le C_{\fcb} &\sum_{r=1}^{k-1}(a_{r+1} - b_{r})^{-\gamma}\Big(|\varphi|_\infty+\Lip(\varphi)\sum_{i=1}^{k}\sum_{j=1}^{b_i - a_i + 1}\Lip(\Phi_i;j)\Big).
\end{align*}
\end{prop}
Note that here $[-R,R]^k$ is endowed with the $\ell^1$ metric, i.e. $|(x_1,\dots,x_k)|=\sum_{i=1}^k |x_i|$.

\begin{proof}
By the argument in the proof of~\cite[Lemma~4.1]{F22}, this follows from Definition~\ref{de:fcb}.
\end{proof}

\subsection{The main result}

Let $v_n:M\to\R$ be a sequence of Lipschitz functions such that $\sup_n\|v_n\|_{\Lip}<\infty$, where $\|v\|_{\Lip}=|v|_\infty+\Lip(v)$.
Denote $S_{j,n}:=\sum_{i=j}^{n-1}\bar v_i\circ \cal T^i$, $\Sigma_{j,n}^2:=\E_m[(\sum_{i=j}^{n-1}\bar v_i\circ \cal T^i)^2]$, where
$\bar v_i:=v_i-\int_M v_i\circ \cal T^i\rmd m$.
In particular, denote $S_n=\sum_{i=0}^{n-1}\bar v_i\circ \cal T^i$, $\Sigma_n^2=\E_m[(\sum_{i=0}^{n-1}\bar v_i\circ \cal T^i)^2]$.
For $t\in[0,1]$, set
\[
N_n(t):=\max\{0 \le j\le n: \Sigma_j^2 \le t\Sigma_n^2\}.
\]
Let $D[0,1]$ denote the space of c\`adl\`ag functions, which we endow with the metric induced by the supremum norm. 
Let $n\ge 1$ and assume that $\Sigma_n>0$. 
Consider the process $W_{n}\in D[0,1]$ defined by
\begin{equation}\label{wnt}
W_{n}(t):=\frac{1}{\Sigma_n}\sum_{i=0}^{N_n(t)-1}\bar v_i\circ \cal T^i,\quad t\in[0,1].
\end{equation}

Now, we introduce our main result on convergence rates in the weak invariance principle with respect to the Wasserstein-$r$ metric.
	\begin{thm}\label{thnon}
		Suppose that $\{\cal T^n\}_{n\ge 0}$ satisfies the Functional Correlation Bound with rate $n^{-\gamma}$. 
		Let $\{v_n\}_{n\ge 0}$ be a sequence of Lipschitz functions with $\sup_n\|v_n\|_{\Lip}<\infty$.     
        Let $\gamma>1$, $1\le r<\gamma+1$ and $\eps>0$. 
        Then there exists a constant $C>0$ such that
		\begin{align*}
            \cal W_r(W_n,B)\le\begin{cases}
                 Cn^\eps \big[(n\Sigma_n^{-(1+2\gamma)})^{\ell(\gamma,r)}+\Sigma_n^{-\frac{\gamma-1}{2\gamma}}n^{(1-\frac{\gamma}{2})\frac{\gamma-1}{2\gamma}}+ n^{\frac{1}{\gamma+1}}\Sigma_n^{-1}\big], & 1<\gamma<2,\\
                Cn^\eps \big[(n\Sigma_n^{-(1+2\gamma)})^{\ell(\gamma,r)} + n^{\frac{1}{\gamma+1}}\Sigma_n^{-1} \big], & \gamma\ge 2,
            \end{cases}
		\end{align*}
		where $B$ is a standard Brownian motion and 
		\begin{align*}
			\ell(\gamma,r)= \frac{(\gamma-1)(\gamma+1-r)}{4(\gamma+1)(\gamma^2+\gamma+1)-r(3\gamma^2+5\gamma+4)}.
		\end{align*}
        In particular, if $\lim_{n\to\infty}\Sigma_n/n^\kappa=\infty$ for some $\kappa>1/(\gamma+1)$, then $\cal W_r(W_n,B)\to 0$ as $n\to \infty$ and thus $W_n$ weakly converges to $B$.
	\end{thm}

\begin{rem}
(i) In the variance growth assumption $\lim_{n\to\infty}\Sigma_n/n^\kappa{=\infty}$, the hypothesis that $\kappa>1/(\gamma+1)$ is essentially sharp. In the appendix, we provide an example to show that the CLT (and hence the WIP) can fail for $\kappa<1/(\gamma+1)$.

(ii) Note that $\lim_{\gamma\to \infty}\ell(\gamma,r)=0$ and
$\lim_{\gamma\to \infty}(1+2\gamma)\ell(\gamma,r)=\tfrac12$.
Hence if the Functional Correlation Bound holds with a superpolynomial rate and 
the variance $\Sigma_n^2$ grows polynomially, then Theorem~\ref{thnon} yields that $\cal W_r(W_n,B)=O(\Sigma_n^{-(1/2-\eps)})$ for all $\eps>0$, $1\le r<\infty$.
In~\cite{LW24}, the same convergence rate was also obtained for sequential expanding maps without any variance growth assumption.
Even in the stationary case, obtaining a better rate than $O(\Sigma_n^{-1/2})$ in the Wasserstein metric is an open problem for dynamical systems.

(iii) Theorem~\ref{thnon} also yields convergence rates in the L\'{e}vy-Prokhorov distance $\pi$. 
Indeed, for any probability measures $\mu$ and $\nu$, we have $\pi(\mu,\nu)\le \mathcal{W}_{p}(\mu,\nu)^{\frac{p}{p+1}}$ for $1\le p<\infty$.
\end{rem}    

\begin{rem}
Note that $\cal W_r$ defines a metric on the space of Radon probability measures on $D[0,1]$ with finite $r$-th moment. 
The law of $W_n$ is Radon since $W_n = \rho(\bar v_0,\dots,\bar v_{n-1}\circ \cal T^{n-1})$ where $\rho:\R^n\to D[0,1]$ is continuous. 
Likewise, the law of $B$ is Radon since Brownian motion has continuous sample paths. 
By similar reasoning, all random elements of $D[0,1]$ considered in this paper have Radon laws.
\end{rem}

\section{Examples}\label{sec:examples} In this section, we give examples of sequential dynamical systems to which our main result applies.

\subsection{Sequential intermittent maps}

Let $M=[0,1]$. We consider the family of Liverani--Saussol--Vaienti (LSV) maps $T_\alpha:M\rightarrow M$ introduced in~\cite{LSV99}, given by
\begin{align}\label{P-M}
  T_\alpha(x)=
  \begin{cases}
  x(1+2^\alpha x^\alpha) & x\in[0,\frac{1}{2}),\\
  2x-1  & x\in[\frac{1}{2},1],
  \end{cases}
\end{align}
where $\alpha\in (0,1)$ is a parameter. 
It follows from \cite{LSV99} that each $T_\alpha$ admits a unique absolutely continuous invariant probability measure $\mu_\alpha$ whose density belongs to the convex cone
\begin{align*}
\cal C_*(\alpha)=\{f\in C((0,1])\cap L^1(\Leb): & f\ge 0, f\text{~decreasing},\\
                                        & x^{\alpha+1}f \text{~increasing}, f(x)\le 2^{\alpha}(2+\alpha)x^{-\alpha}\Leb(f)\}.
\end{align*}
Note that $1\in \cal C_*(\alpha)$ and $\cal C_*(\alpha)\subset \cal C_*(\beta)$ for any $0<\alpha<\beta<1$. 
For $n\ge 1$, let $T_n=T_{\alpha_n}$ where $\alpha_n\in (0,1)$ with $\alpha_*=\sup_n \alpha_n<1$. 
Let $m$ be an absolutely continuous probability measure with density in $\cal C_*(\alpha_*)$. 
Then, by~\cite{L17}, the Functional Correlation Bound is satisfied with rate $n^{-\gamma}$ for all $\gamma<1/\alpha_*-1$ (more precisely, with rate $n^{1-1/\alpha_*}(\log n)^{1/\alpha_*}$) so Theorem~\ref{thnon} yields rates in the WIP whenever $\alpha_*<\frac12$ and $\Sigma_n\gtrsim n^\kappa$ for some $\kappa>\alpha_*$. 
The variance growth assumption is weaker than that in previous results,
which in particular require $\kappa>\frac14$ for the CLT~\cite{NTV18} and $\kappa>\frac13$ for the WIP~\cite{leppanen_nakajima_nakano}. 
Moreover, the restriction $\alpha_*<\frac12$ is essentially sharp. Indeed, even in the stationary setting where $v_n\equiv v$, $T_n\equiv T_{\alpha}$, the CLT fails whenever $\alpha>\frac12$ and $v(0)\ne \int_M v\, \rmd \mu_\alpha$~\cite{gouezel_stable}.

\subsection{Random intermittent maps}\label{subsec:random_lsv}
In particular, we obtain rates in the quenched WIP for random LSV maps over an ergodic base.
Let $0<\alpha_*<\frac12$ and define $\Omega=(0,\alpha_*]^{\Z}$. Let $\sigma:\Omega\to \Omega$ be the shift map, $(\sigma\omega)_i = \omega_{i+1}$. 
Let $\P$ be an ergodic invariant probability measure on $\Omega$. Define $T_\omega=T_{\omega_0}$ for $\omega\in \Omega$. 
By the proof of~\cite[Corollary 3.3]{Su22}, there exists a family of absolutely continuous probability measures $\mu_\omega$ such that $(T_\omega)_* \mu_\omega= \mu_{\sigma \omega}$ and $ \rmd\mu_\omega/\rmd\Leb\in \cal C_*(\alpha_*)$ for a.e.\ $\omega\in \Omega$. 
Set $T_\omega^n = T_{\sigma^{n-1}\omega}\circ \dots\circ T_\omega$.

Let $v:\Omega\times M\to \R$ be a measurable function such that $v_\omega:=v(\omega,\cdot)$ is Lipschitz for all $\omega\in \Omega$, $\int_M v_\omega \rmd\mu_\omega=0$ and $\sup_\omega\|v_\omega\|_{\Lip}<\infty$. 
For $\omega\in \Omega$, let $\Sigma_n^2(\omega)=\E_{\mu_\omega}[(\sum_{i=0}^{n-1}v_{\sigma^i \omega}\circ T_\omega^i)^2]$ and \[N_n^\omega(t)=\max\{0\le j\le n: \Sigma_j^2(\omega)\le t\Sigma_n^2(\omega)\}.\]
Given $n\ge 1$ such that $\Sigma_n(\omega)>0$, define $W_n^\omega\in D[0,1]$ by 
\[ W_n^\omega(t)=\frac{1}{\Sigma_n(\omega)}\sum_{i=0}^{N_n^\omega(t)-1}v_{\sigma^i \omega}\circ T^i_\omega.\]
Let $\mu$ be the probability measure on $\Omega\times M$ defined by $\rmd\mu(\omega,x)=\rmd\mu_\omega(x)\, \rmd\P(\omega)$. 
Recall that $v$ is called a coboundary if there exists $\chi \in L^1(\Omega\times M,\mu)$ such that $v_\omega = \chi_{\sigma \omega}\circ T_\omega -\chi_\omega$ holds $\mu_\omega$-a.s.\ for a.e.\ $\omega$. 
\begin{cor}
    Suppose that $v$ is not a coboundary. Set $\gamma_*=1/\alpha_*-1$. Let $1\le r<\gamma_*+1$ and $\eps>0$. Then for $\P$-a.e.\ $\omega \in \Omega$, there exist constants $C_\omega>0$ and $n_0(\omega)$ such that for $n\ge n_0(\omega)$,
    \begin{equation*}
        \cal W_r(W_n^\omega, B)\le \begin{cases}
            C_\omega n^\eps\big[n^{-(\gamma_*-\frac12)\ell(\gamma_*,r)}+n^{-\frac{(\gamma_*-1)^2}{4\gamma_*}}\big], & 1<\gamma_* \le 2,\\
            C_\omega n^{\eps - (\gamma_*-\frac12)\ell(\gamma_*,r)}, & \gamma_* > 2.
        \end{cases}
    \end{equation*}
\end{cor}
\begin{proof}
Let $L_a$ denote the transfer operator corresponding to $T_a$. By the proof of~\cite[Corollary~3.3]{Su22}, $\rmd\mu_\omega/\rmd\Leb$ is the limit of $L_{\omega_{-1}}\dots L_{\omega_{-n}}1$ in $L^1(\Leb)$ for a.e. $\omega$. 
Since $\cal C_*(\alpha_*)$ is closed in $L^1$ and $L_\beta \cal C_*(\alpha_*)\subset \cal C_*(\alpha_*)$ for any $0<\beta\le \alpha_*$, it follows that $\rmd\mu_\omega/\rmd\Leb\in \cal C_*(\alpha_*)$.
It follows from the proof of~\cite[Corollary 3.3]{Su22} that the assumptions of~\cite[Theorem 2.10]{Su22} hold. 
Consequently, by the argument at the end of~\cite[p.4837]{Su22} we have $\lim_{n\to\infty}\Sigma_n^2(\omega)/n>0$ so there exists $n_0(\omega)$ such that $\Sigma^2_n(\omega)\gtrsim n$ for $n\ge n_0(\omega)$.

Set $m=\mu_\omega$, $T_n = T_{\omega_{n-1}}$ and $v_n = v_{\sigma^n \omega}$. Then as explained above, the Functional Correlation Bound is satisfied with rate $n^{-\gamma}$ for any $1<\gamma<\gamma_*$. 
Applying Theorem~\ref{thnon} and using $\Sigma_n^2(\omega)\gtrsim n$ for $n\ge n_0(\omega)$, we obtain the stated bound with $\gamma$ in place of $\gamma_*$, since a straightforward calculation shows that $(\gamma-\frac12)\ell(\gamma,r)\le \frac12-\frac{1}{\gamma+1}$. 
Since $\eps>0$ is arbitrary, the result follows by taking $\gamma$ sufficiently close to $\gamma_*$.
\end{proof}

In the case of iid random maps,~\cite{liu_saussol_vaienti_wang} also obtained rates in the quenched WIP in the Wasserstein-$p$ metric. More precisely, their result covers the case where $\P = \nu^{\otimes \Z}$ and $\nu$ is the uniform measure on $[\alpha_0,\alpha_*]$ with $0<\alpha_0<\frac15$ and $\alpha_*<1$.

\subsection{Sequential heterochaos baker maps} Heterochaos baker maps were introduced in~\cite{saiki_sanjuan_yorke_18,saiki_takahasi_yorke_21} as a simple example with partially hyperbolic dynamics.
     Let $M=[0,1]^2$ and fix an integer $N\ge 2$. For $a\in (0,\frac1N)$, define a piecewise linear expanding map $\tau_a:[0,1]\to [0,1]$ by 
    \begin{equation*} \tau_a(x) = \begin{cases}
        \dfrac{x-(i-1)a}{a}, & x \in [(i-1)a,ia),\ i = 1,\dots,N, \\[1.5ex]
        \dfrac{x-Na}{1-Na}, & x\in [aN, 1].
    \end{cases}\end{equation*}
    Let $\{z\}=z-[z]$ denote the fractional part of $z\in \R$. Define $g_a:M \to [0,1]$ by 
    \begin{equation*}
        g_a(x,y)=\begin{cases}
            \dfrac{y}{N}+\dfrac{i-1}{N}, & x \in [(i-1)a,ia), \ i = 1,\dots,N,\\
            \{Ny\}, & x\in [aN, 1],\  y<1,\\
            1, & x \in [aN,1],\ y=1.
        \end{cases}
    \end{equation*}
    The heterochaos baker map $T_a:M\to M$ is defined by $T_a(x,y)=(\tau_a(x),g_a(x,y))$. The map $T_a$ preserves Lebesgue measure and is ergodic. It has expanding center (i.e., the central Lyapunov exponent is positive) if and only if $a\in (0,\frac{1}{2N})$. For $n\ge 1$, let $a_n\in (0,\frac{1}{2N})$ and define $T_n = T_{a_n}$. Suppose that $\sup_n a_n <\frac{1}{2N}$ and 
    \begin{equation}\label{eq:hetero_admissible}
        \max\bigg\{0, \frac{4N \inf_n a_n - 1}{4N^2\sup_n a_n} \bigg\}<\inf_n a_n.
    \end{equation}
    Note that~\eqref{eq:hetero_admissible} is satisfied if $\inf_n a_n\in (0,\frac{1}{4N})$, or if $a_n$ is sufficiently close to a constant $a\in (0,\frac{1}{2N})$ for all $n$.
    Let $m$ be a probability measure with density that is Lipschitz and bounded away from zero. Then the Functional Correlation Bound is satisfied with a stretched exponential rate by~\cite[Theorem 1.2]{leppanen_takahasi} so Theorem~\ref{thnon} applies and $\gamma$ can be taken arbitrarily large. In particular, if $\Sigma_n\gtrsim n^\kappa$ for some $\kappa>0$, then we obtain that $\cal W_r(W_n,B)\lesssim \Sigma_n^{-(1/2-\eps)}$ for any $1\le r<\infty$ and $\eps>0$. 

\subsection{Further examples} By~\cite{F22}, if $T_n \equiv T$ is a mixing nonuniformly hyperbolic map modeled by a Young tower with $O(n^{-(\gamma+1)})$ tails, $\gamma>0$, then the Functional Correlation Bound holds with rate $n^{-\gamma}$. In particular, Theorem~\ref{thnon} applies to examples such as Anosov diffeomorphisms, the Sinai billiard map and Bunimovich flowers. In~\cite{F25}, Wasserstein rates were obtained in this setting in the stationary case where $v_n\equiv v$. 

By~\cite[Proposition 3.5]{leppanen_stenlund-sunklodas}, the Functional Correlation Bound is satisfied with an exponential rate for random Lasota--Yorke maps. Hence, similarly to Subsection~\ref{subsec:random_lsv}, we obtain rates in the quenched WIP whenever the variance grows linearly (which holds under a non-coboundary condition; see~\cite[Theorem 3.3]{leppanen_stenlund-sunklodas}).

\section{Moment estimates}\label{sec:moment_estimates}
This section is dedicated to the proof of the following moment bound, which will be key in the proof of our main theorem.
\begin{prop}\label{prop:moment_bds_fcb}
    Suppose that $\{\cal T^n\}_{n\ge 0}$ satisfies the Functional Correlation Bound with rate $n^{-\gamma}$, $\gamma>1$ and constant $C_{\fcb}>0$. 
	Let $\{v_n\}_{n\ge 0}$ be a sequence of Lipschitz functions with $\sup_n\|v_n\|_{\Lip}<\infty$.
    Let $\eps>0$. Then there exists a constant $D>0$ such that 
    \begin{equation}\label{eq:induc_moment_claim}
    	\Big\|\max_{1\le j\le n}|S_{m,m+j}| \Big\|_{\gamma+1}^{\gamma+1} \le Dn^\eps\big((\Var S_{m,m+n})^{(\gamma+1)/2}+n\big)
    \end{equation}
    for all $m\ge 0$, $n\ge 1$.
    Moreover, the constant $D$ depends only on $C_{\fcb}$, $\gamma$, $\sup_n \|v_n\|_{\Lip}$ and $\eps$.

\end{prop}
Throughout the remainder of this section, we fix  
$\{\cal T^n\}_{n\ge 0}$ that satisfies the Functional Correlation Bound with rate $n^{-\gamma}$, $\gamma>1$ and constant $C_{\fcb}>0$, and a sequence of Lipschitz functions $v_n:M\to \R$ such that $\sup_n \|v_n\|_{\mathrm{Lip}}<\infty$. In this section, we denote by $C$ various constants that depend only on $C_{\fcb}$, $\gamma$, $\sup_n \|v_n\|_{\Lip}$ and $\eps$; the same convention applies to the implicit constant in $\lesssim$.

We prove Proposition~\ref{prop:moment_bds_fcb} by induction on $n$; the proof is inspired by the proof of \cite[Theorem 2.4(a)]{F22}. Before completing the proof we need several preparatory lemmas; the following lemma will also be useful elsewhere.

\begin{lem}\label{lem:covar_approx_block_sum}
		For all $n\ge 1$, $k\ge 1$ and $m=\ell_0<\ell_1<\dots<\ell_k=n$, we have 
		\[ \bigg|\Sigma^2_{m,n}-\sum_{i=0}^{k-1}\Sigma^2_{\ell_i,\ell_{i+1}}\bigg|\le C\sum_{i=0}^{k-1}\big((\ell_{i+1}-\ell_i)^{2-\gamma}+\log(\ell_{i+1}-\ell_i)\big). \]
    \end{lem}
	\begin{proof}
		First note that
		\begin{align*}
			\Sigma_{m,n}^2 = \sum_{i=0}^{k-1}\E_m[S_{\ell_i,\ell_{i+1}}\, S_{m,n} ]=\sum_{i=0}^{k-1}\Big(\Sigma_{\ell_i,\ell_{i+1}}^2 + \E_m\big[S_{\ell_i,\ell_{i+1}}(S_{m,\ell_i}+S_{\ell_{i+1},n})\big]\Big).
		\end{align*}
        Now, for any $j_1<j_2<j_3$, we have
        \begin{align*}
            \big|\E_m[S_{j_1,j_2} \, S_{j_2,j_3} ]\big|&=\bigg|\sum_{a=j_1}^{j_2-1}\sum_{b=j_2}^{j_3-1}\E_m[\bar v_a\circ\cal{T}^a\, \bar v_b\circ\cal{T}^b ]\bigg|.
        \end{align*}
        Observe that \[\E_m[\bar v_a\circ\cal{T}^a\, \bar v_b\circ\cal{T}^b ]=\int v_a \circ \cal T^a v_b \circ \cal T^b \rmd m-\int v_a\circ \cal T^a \rmd m\int v_b \circ \cal T^b \rmd m,\]
        so by applying the Functional Correlation Bound with $G(x_1,x_2)=v_a(x_1)v_b(x_2)$ we obtain that $|\E_m[\bar v_a\circ\cal{T}^a\, \bar v_b\circ\cal{T}^b ]|\lesssim (b-a)^{-\gamma}.$
        It follows that 
        \begin{align*}
            \big|\E_m[S_{j_1,j_2} \, S_{j_2,j_3} ]\big|             &\lesssim \sum_{a=j_1}^{j_2-1}\sum_{b=j_2}^{j_3-1}(b-a)^{-\gamma}
\lesssim \sum_{a=j_1}^{j_2    
            -1}(j_2 - a)^{1-\gamma}\lesssim (j_2-j_1)^{2-\gamma}+\log(j_2-j_1).
        \end{align*}
        A similar argument shows that we also have
        \begin{equation}\label{eq:sums_corr_bd}
            \big|\E_m[S_{j_1,j_2} \, S_{j_2,j_3}]\big|\lesssim (j_3-j_2)^{2-\gamma}+\log(j_3-j_2).
        \end{equation}
        Hence the result follows by using the bound
        \begin{align*}
            \Big|\Sigma_{m,n}^2 - \sum_{i=0}^{k-1}\Sigma_{\ell_i,\ell_{i+1}}^2\Big|&\le \sum_{i=0}^{k-1}\big|\E_m[S_{\ell_i,\ell_{i+1}}(S_{m,\ell_i}+S_{\ell_{i+1},n})]\big|.
        \end{align*}
	\end{proof}
    The following result will be useful in Section~\ref{sec:main_thm_pf}.
    \begin{lem}\label{lem:var_diff}
        We have $|\Sigma_n^2 - \Sigma_m^2|\le C(n-m)$ for all $n> m\ge 0$.
    \end{lem}
    \begin{proof}
        Since $|\E_m[\bar v_a\circ \cal T^a\, \bar v_b\circ \cal T^b]|\lesssim (b-a)^{-\gamma}$ and $\gamma>1$, a standard argument shows that $\Sigma_{m,n}^2\lesssim n-m$. By~\eqref{eq:sums_corr_bd}, it follows that 
        \[ |\Sigma_n^2 - \Sigma_m^2|=|\Sigma_{m,n}^2+2\E_m[S_m S_{m,n}]|\lesssim n-m +(n-m)^{2-\gamma}+\log(n-m)\lesssim n-m. \qedhere \]
    \end{proof}
    Let $k\ge 1$ and $N\ge 2k$. Define $\alpha(i) = m+\big[\tfrac{iN}{2k}\big]$ for $i\ge 0$. Note that 
    \begin{equation}\label{eq:a_i-bounds}
    	\tfrac{N}{2k}-1\le \alpha(i+1) - \alpha(i) \le \tfrac{N}{2k}+1\le \tfrac{N}{k}.
    \end{equation}
    Let $X_i = S_{\alpha(2i),\alpha(2i+1)}$. Let $\widehat X_0,\dots,\widehat X_{k-1}$ be independent random variables such that $\widehat X_i =_d X_i $ for $0\le i\le k-1$.

\begin{lem}\label{lem:moments_indep_approx} 
	For all $p\in [1,\infty)$ and integers $N\ge 2k$, $k\ge 1$, we have
    \[ 
    \E_m\bigg[\max_{\ell\le k}\bigg|\sum_{i=0}^{\ell-1}X_i\bigg|^p\bigg]\le C^p pk^{\gamma+1}N^{p-\gamma}+ \E\bigg[\max_{\ell\le k}\bigg|\sum_{i=0}^{\ell-1}\widehat X_i\bigg|^p\bigg].
    \]
\end{lem}
\begin{proof}
	This follows from similar calculations to the proof of~\cite[Lemma 5.4]{F22}; we give the details for completeness.
    Note that $X_i(x)=\Phi_i(\cal T^{a_i}x,\ldots,\cal T^{b_i}x)$, where $a_i = \alpha(2i)$, $b_i = \alpha(2i+1)-1$ and
   \[
     \Phi_i(y_{a_i},\ldots,y_{b_i})=\sum_{j=a_i}^{b_i}\bar v_j(y_j).
    \]
   Define $R=\max_i|\Phi_i|_{\infty}$. Then 
   \begin{align*}
       \Delta = \E_m\bigg[\max_{\ell\le k}\bigg|\sum_{i=0}^{\ell-1}X_i\bigg|^p\bigg] - \E\bigg[\max_{\ell\le k}\bigg|\sum_{i=0}^{\ell-1}\widehat X_i\bigg|^p\bigg]=\E_m[\varphi(X_0,\dots,X_{k-1})]-\E[\varphi(\widehat X_0,\dots,\widehat X_{k-1})],
   \end{align*}
   where $\varphi:[-R,R]^k\to \R$ is given by $\varphi(y_0,\ldots,y_{k-1})=\max_{\ell\le k}\big|\sum_{i=0}^{\ell-1}y_i\big|^p$. 
We proceed by using Proposition~\ref{lip} to bound $|\Delta|$.
    First note that $\Lip(\Phi_i;j) = \Lip\, \bar v_{a_i + j - 1} =\Lip\, v_{a_i + j - 1}$ and \[|\Phi_i|_\infty\le (b_i - a_i + 1)\sup_n |\bar v_n|_\infty\le 2(b_i - a_i +1)\sup_n |v_n|_\infty.\] By~\eqref{eq:a_i-bounds}, $b_i - a_i +1= \alpha(2i+1)-\alpha(2i)\le \frac{N}{k}$. Let $V=\sup_n \|v_n\|_{\Lip}$. It follows that $R\le 2VN/k$ and 
    \begin{equation*}
        \sum_{j=1}^{b_i - a_i + 1}\Lip(\Phi_i;j)\le VN/k.
    \end{equation*}
   Write $\varphi=\max\{h_1,\dots,h_k\}$ where $h_\ell(y_0,\dots,y_{k-1})=\big|\sum_{i=0}^{\ell-1}y_i\big|^p$. By~\cite[Proposition 5.3]{F22}, $|h_\ell|_\infty\le (kR)^p$ and $\Lip(h_\ell)\le p(kR)^{p-1}$ and so $|\varphi|_\infty\le (kR)^p$. Moreover, since 
   \[ |\varphi(\bar y)-\varphi(\bar y')|\le \max_{\ell \le k}|h_\ell(\bar y)-h_\ell(\bar y')|\]
   for all $\bar y,\bar y'\in [-R,R]^k$, it follows that $\Lip(\varphi)\le p(kR)^{p-1}$. Finally, note that $a_i - b_{i-1} =\alpha(2i)-\alpha(2i-1)+1\ge \frac{N}{2k}$ by~\eqref{eq:a_i-bounds}. Hence by applying Proposition~\ref{lip} we obtain
    \begin{align*}
        |\Delta|&\lesssim \sum_{r=1}^{k-1}(a_r - b_{r-1})^{-\gamma}\bigg(|\varphi|_\infty+\Lip(\varphi)\sum_{i=0}^{k-1}\sum_{j=1}^{b_i - a_i + 1}\Lip(\Phi_i;j)\bigg)\\
        &\lesssim k\big(\tfrac{N}{2k}\big)^{-\gamma}\big((kR)^p + p(kR)^{p-1}k\cdot \tfrac{VN}{k}\big)\\
        &\le 2^\gamma k^{1+\gamma}N^{-\gamma}\big((2VN)^p+p(2VN)^{p-1}VN\big)
        \lesssim p (2V)^p k^{\gamma+1}N^{p-\gamma},
    \end{align*}
   as required.
\end{proof}

We also need the following result, which is an immediate consequence of Rosenthal's inequality~\cite[Theorem 3]{rosenthal} and Doob's maximal inequality.
\begin{prop}\label{prop:rosenthal}
Fix  $p\ge 2$. Then there exists a constant $M>0$ such that for all $k\ge1$, for any independent, mean zero random variables $\widehat X_0,\ldots, \widehat X_{k-1}\in L^p$, 
\[
\E\bigg[\max_{\ell \le k}\Big|\sum_{i=0}^{\ell-1}\widehat{X}_i\Big|^p\bigg] \le M\bigg(\Big(\sum_{i=0}^{k-1}\E[\widehat{X}_i^2]\Big)^{p/2} + \sum_{i=0}^{k-1}\E[|\widehat{X}_i|^p]\bigg).
\]
\end{prop}
The following lemma is the main step in the proof of Proposition~\ref{prop:moment_bds_fcb}.
\begin{lem}\label{lem:moments_induc_claim}
    Let $D>0$, $k\ge 1$ and $N\ge N_0(k,\gamma)$ be such that~\eqref{eq:induc_moment_claim} holds for all $m\ge 0$ and $n < N$. Then 
    \[ \Big\|\max_{1\le j\le N}|S_{m,m+j}| \Big\|_{\gamma+1}^{\gamma+1}\le C(k^{\gamma+1}+Dk^{-\eps})N^\eps\big((\Var S_{m,m+N})^{(\gamma+1)/2}+N\big)\]
    for all $m\ge 0$. Here $N_0(k,\gamma)$ depends only on $k$ and $\gamma$.
\end{lem}
\begin{proof}
	First note that 
	\begin{align*}
		\max_{1\le j\le N}|S_{m,m+j}|&\le \max_{\ell\le 2k}\bigg|\sum_{i=0}^{\ell-1}S_{\alpha(i),\alpha(i+1)}\bigg|+\max_{\ell\le 2k-1}\max_{\alpha(\ell)\le j\le \alpha(\ell+1)}|S_{\alpha(\ell),j}|\\
		&\le I_1 + I_2+M,
	\end{align*}
	where 
	\begin{align*}
		I_1 &= \max_{\ell\le k}\bigg|\sum_{i=0}^{\ell - 1}S_{\alpha(2i),\alpha(2i+1)}\bigg| ,\qquad I_2 = \max_{\ell \le k} \bigg|\sum_{i=0}^{\ell-1}S_{\alpha(2i+1),\alpha(2i+2)}\bigg|,\\
        M &= \max_{\ell\le 2k-1}\max_{\alpha(\ell)\le j\le \alpha(\ell+1)}|S_{\alpha(\ell),j}|.
	\end{align*}
	Now 
	\begin{align}\label{eq:M_bound}
		\E_m[M^{\gamma+1}]\le \sum_{\ell=0}^{2k-1}\E_m \bigg[\max_{\alpha(\ell)\le j\le \alpha(\ell+1)}|S_{\alpha(\ell),j}|^{\gamma+1}\bigg].
	\end{align}
	
    Let $N\ge 2k$. Then, by applying Lemma~\ref{lem:moments_indep_approx} with $p=\gamma+1$ we obtain that 
    \[ \E_m[I_1^{\gamma+1}]\lesssim k^{\gamma+1}N + \E\bigg[\max_{\ell\le k}\bigg|\sum_{i=0}^{\ell-1}\widehat X_i\bigg|^{\gamma+1}\bigg]. \]
    By Proposition~\ref{prop:rosenthal},
    \begin{align*}
    	\E\bigg[\max_{\ell\le k}\bigg|\sum_{i=0}^{\ell-1}\widehat X_i\bigg|^{\gamma+1}\bigg] &\lesssim \Big(\sum_{i=0}^{k-1}\E[\widehat{X}_i^2]\Big)^{\frac{\gamma+1}{2}} + \sum_{i=0}^{k-1}\E[|\widehat{X}_i|^{\gamma+1}]\\
    	&=\Big(\sum_{i=0}^{k-1}\Sigma_{\alpha(2i),\alpha(2i+1)}^2\Big)^{\frac{\gamma+1}{2}} + \sum_{i=0}^{k-1}\E_m[|S_{\alpha(2i),\alpha(2i+1)}|^{\gamma+1}]
    \end{align*}
    and so we obtain that 
    \[ \E_m[I_1^{\gamma+1}] \lesssim k^{\gamma+1}N + \Big(\sum_{i=0}^{k-1}\Sigma_{\alpha(2i),\alpha(2i+1)}^2\Big)^{\frac{\gamma+1}{2}} + \sum_{i=0}^{k-1}\E_m[|S_{\alpha(2i),\alpha(2i+1)}|^{\gamma+1}]. \]
    By the same argument,
    \[ \E_m[I_2^{\gamma+1}] \lesssim k^{\gamma+1}N + \Big(\sum_{i=0}^{k-1}\Sigma_{\alpha(2i+1),\alpha(2i+2)}^2\Big)^{\frac{\gamma+1}{2}}+\sum_{i=0}^{k-1}\E_m[|S_{\alpha(2i+1),\alpha(2i+2)}|^{\gamma+1}]. \]
    Combining these bounds yields that
    \begin{align*}
    	\E_m[I_1^{\gamma+1}]+\E_m[I_2^{\gamma+1}]\lesssim k^{\gamma+1}N+ \Big(\sum_{i=0}^{2k-1}\Sigma_{\alpha(i),\alpha(i+1)}^2\Big)^{\frac{\gamma+1}{2}}+\sum_{i=0}^{2k-1}\E_m[|S_{\alpha(i),\alpha(i+1)}|^{\gamma+1}].
    \end{align*}
    By~\eqref{eq:M_bound}, it follows that 
    \begin{align}
    	\E_m\Big[\max_{1\le j\le N}|S_{m,m+j}|^{\gamma+1}\Big]&\le\E_m[|M+I_1+I_2|^{\gamma+1}] \le 3^{\gamma+1}\big(\E_m[M^{\gamma+1}]+\E_m[I_1^{\gamma+1}]+\E_m[I_2^{\gamma+1}]\big) \nonumber \\
    	&\lesssim k^{\gamma+1}N + \Big(\sum_{i=0}^{2k-1}\Sigma_{\alpha(i),\alpha(i+1)}^2\Big)^{\frac{\gamma+1}{2}}+\sum_{\ell=0}^{2k-1}\E_m\Big[\max_{\alpha(\ell)\le j\le \alpha(\ell+1)}|S_{\alpha(\ell),j}|^{\gamma+1}\Big].   \label{eq:M_Ij_sum}	
    \end{align}
    
    	Now, by the inductive hypothesis~\eqref{eq:induc_moment_claim}, we have 
    \begin{align*}
    	\sum_{\ell=0}^{2k-1}\E_m\Big[\max_{\alpha(\ell)\le j\le \alpha(\ell+1)}|S_{\alpha(\ell),j}|^{\gamma+1}\Big]&\le D\sum_{\ell=0}^{2k-1}(\alpha(\ell+1)-\alpha(\ell))^\eps\big(\Sigma_{\alpha(\ell),\alpha(\ell+1)}^{\gamma+1}+\alpha(\ell+1)-\alpha(\ell)\big)\\
    	&\le D\big(\tfrac{N}{k}\big)^\eps \bigg( \sum_{\ell=0}^{2k-1}\Sigma_{\alpha(\ell),\alpha(\ell+1)}^{\gamma+1}+N \bigg),
    \end{align*}
    where we have used that $\alpha(\ell+1)-\alpha(\ell)\le N/k$ by~\eqref{eq:a_i-bounds} and $\sum_{\ell=0}^{2k-1}(\alpha(\ell+1)-\alpha(\ell))=N$. Since $(\gamma+1)/2>1$, it follows that
    \begin{equation*}
    	\sum_{\ell=0}^{2k-1}\E_m\Big[\max_{\alpha(\ell)\le j\le \alpha(\ell+1)}|S_{\alpha(\ell),j}|^{\gamma+1}\Big]\le D\big(\tfrac{N}{k}\big)^\eps \bigg( \bigg(\sum_{\ell=0}^{2k-1}\Sigma_{\alpha(\ell),\alpha(\ell+1)}^2\bigg)^{\frac{\gamma+1}{2}}+N \bigg).
    \end{equation*}
    Hence by combining this bound with~\eqref{eq:M_Ij_sum}, we obtain that
    	\begin{align*}
    	\E_m\Big[\max_{1\le j\le N}|S_{m,m+j}|^{\gamma+1}\Big]\lesssim k^{\gamma+1}N + \big(1+D\big(\tfrac{N}{k}\big)^\eps\big)\bigg( \bigg(\sum_{\ell=0}^{2k-1}\Sigma_{\alpha(\ell),\alpha(\ell+1)}^2\bigg)^{\frac{\gamma+1}{2}}+N \bigg).
    	\end{align*}
        
        Finally, by Lemma~\ref{lem:covar_approx_block_sum} and again using that $\alpha(\ell+1)-\alpha(\ell)\le N/k$, we have
    \begin{align*}
        \sum_{i=0}^{2k-1}\Sigma_{\alpha(i),\alpha(i+1)}^2&\lesssim \Sigma_{m,m+N}^2+\sum_{i=0}^{2k-1}\big((N/k)^{2-\gamma}+\log(N/k)\big)\\
        &\lesssim \Sigma_{m,m+N}^2+N^{2-\gamma}k^{\gamma-1}+k\log N.
    \end{align*}
    Now, since $\frac{1}{2}(2-\gamma)(\gamma+1)<1$ for all $\gamma>1$, there exists $N_0(k,\gamma)$ such that 
    \[(N^{2-\gamma}k^{\gamma-1}+k\log N)^{\frac{\gamma+1}{2}}\le N\]
    whenever $N\ge N_0(k,\gamma)$.
    It follows that
    \begin{align*}
    		\E_m\Big[\max_{1\le j\le N}|S_{m,m+j}|^{\gamma+1}\Big]
            &\lesssim k^{\gamma+1}N+\big(1+D\big(\tfrac{N}{k}\big)^\eps\big)(\Sigma_{m,m+N}^{\gamma+1}+N)\\
    		&\le (k^{\gamma+1}+1+Dk^{-\eps})N^\eps(\Sigma_{m,m+N}^{\gamma+1}+N),
    \end{align*}
    which concludes the proof.
\end{proof}
We are now in a position to complete the proof of Proposition~\ref{prop:moment_bds_fcb}.
\begin{proof}[Proof of Proposition~\ref{prop:moment_bds_fcb}]We use induction on $n$ to prove that for $D$ large enough,~\eqref{eq:induc_moment_claim} holds for all $n\ge 1$ and $m\ge 0$.
	Let $C>0$ be as in Lemma~\ref{lem:moments_induc_claim} and choose $k$ large enough such that $ Ck^{-\eps}\le \frac12$. Since 
	\[ \Big\|\max_{1\le j\le n}|S_{m,m+j}|\Big\|_{\gamma+1}\le n\sup_j |\bar v_j|_\infty\le 2n \sup_j |v_j|_\infty,\] 
	for $D>0$ large enough~\eqref{eq:induc_moment_claim} holds for all $n\le N_0(k,\gamma)$ and $m\ge 0$. Let $N\ge N_0(k,\gamma)$ and assume that~\eqref{eq:induc_moment_claim} holds for all $n<N$ and $m\ge 0$. Then since $ Ck^{-\eps}\le \frac12$, by Lemma~\ref{lem:moments_induc_claim}, for all $m\ge 0$ we have
	\begin{align*}
		\Big\|\max_{1\le j\le N}|S_{m,m+j}|\Big\|_{\gamma+1}^{\gamma + 1}&\le (Ck^{\gamma+1}+\tfrac12 D)N^\eps((\Var S_{m,m+N})^{(\gamma+1)/2}+N)\\
		&\le DN^\eps((\Var S_{m,m+N})^{(\gamma+1)/2}+N)
	\end{align*}
	whenever $D\ge 2Ck^{\gamma+1}$. Hence by induction, for $D>0$ large enough~\eqref{eq:induc_moment_claim} holds for all $n\ge 1$ and $m\ge 0$.
\end{proof}

\section{Proof of Theorem \ref{thnon}}\label{sec:main_thm_pf}
\setcounter{equation}{0}
In this section, we prove the main result. The proof is based on Bernstein's classical ``big block--small block" method. Recall that for every $t\in[0,1]$, 
\[
N_n(t):=\max\{0 \le j\le n: \Sigma_j^2 \le t\Sigma_n^2\}.
\]
Let $\frac{1}{n}<\eta<\rho\le \frac{1}{4}$. Set $k=[\frac{1}{\rho+\eta}]$ and define $0=t_0<s_1<t_1<\dots<s_k<t_k=1$ by 
\begin{align}
    s_i - t_{i-1}&=\rho, \quad 1\le i\le k,\\
    t_i - s_i &= \eta, \quad 1\le i\le k-1.
\end{align}
Note that $t_k - s_k \le 2(\rho+\eta)$ and since $\rho\le \frac{1}{4}$, we have
\begin{equation}\label{eq:k-asymp}
  \tfrac{1}{4\rho}\le  k\le \tfrac{1}{\rho}\le n.
\end{equation}
For $1\le i\le k$, define big blocks $\cal I_i$ and small blocks $\cal J_i$ by 
\[ \cal I_i = \{N_n(t_{i-1}),\ldots, N_n(s_{i})-1\}, \qquad \cal J_i = \{N_n(s_i),\dots, N_n(t_i)-1\}. \]
Let $\cal I$ denote the union of all big blocks and $\cal J$ denote the union of all small blocks.

Denote
\[
X_i=\frac{1}{\Sigma_n}\sum_{j=N_n(t_{i-1})}^{N_n(s_{i})-1}\bar v_{j}\circ \mathcal T^j, \quad Y_i=\frac{1}{\Sigma_n}\sum_{j=N_n(s_{i})}^{N_n(t_{i})-1}\bar v_{j}\circ \mathcal T^j \quad \text{ for $1\le i\le k$}.
\]
Note that
\begin{equation}\label{eq:blocks_W_n}
    X_i = W_n(s_i)-W_n(t_{i-1}), \qquad Y_i = W_n(t_i)-W_n(s_i)\quad \text{ for $1\le i\le k$}.
\end{equation}
We have defined big blocks and small blocks in this way so that $\E_m[X_i^2]\approx 1/k$ for $1\le i\le k$ (see Proposition~\ref{prop:X_i-var-approx}).

{\bf Step 1. Estimation of the convergence rate between processes $W_n$ and $\widetilde W_n$, where $\widetilde W_n(t)=\sum_{i=1}^{[kt]}X_i$.}

Define
\[
A_n(t)=\frac{1}{\Sigma_n}\sum_{i=0\colon i\in \cal I}^{N_n(t)-1}\bar v_{i}\circ \mathcal T^i.
\]
Then
\[
{\mathcal W}_{\gamma+1}(W_n,\widetilde W_n)\le {\mathcal W}_{\gamma+1}(W_n, A_n)+{\mathcal W}_{\gamma+1}(A_n,\widetilde W_n)=:I+II,
\]
where
\[ I=\cal W_{\gamma+1}(W_n,A_n)\le \frac{1}{\Sigma_n}\bigg\|\max_{m\le n}\Big|\sum_{i=0:i\in \cal J}^{m-1}\bar v_i \circ \cal T^i\Big| \bigg\|_{\gamma+1},\]
\[ II=\cal{W}_{\gamma+1}(A_n,\widetilde W_n)\le 2\big\|\max_{1\le i\le k}|U_i|\big\|_{\gamma+1},\] 
where $U_i =\sup_{t_{i-1}\le r\le s_i}|W_n(r)-W_n(t_{i-1})|$.

\begin{prop}\label{prop:block_var_approx}
		There exists a constant $C>0$ such that \[|\E_m[(W_n(t_2)-W_n(t_1))^2]-(t_2-t_1)|\le C\Sigma_n^{-2}(n^{2-\gamma}+\log n)\] for all $n\ge 1$ and $0\le t_1<t_2\le 1$.
	\end{prop}
	\begin{proof}
		By the definition of $N_n(t)$, we have $\Sigma^2_{N_n(t)}\le t\Sigma^2_n < \Sigma^2_{N_n(t)+1}$ for $t\in [0,1)$ and $\Sigma_{N_n(1)}^2=\Sigma_n^2$. It follows that 
        \[ \Sigma^2_{N_n(t_2)}-\Sigma^2_{N_n(t_1)+1} < (t_2 - t_1)\Sigma^2_n \le \max\{\Sigma^2_{N_n(t_2)+1}-\Sigma^2_{N_n(t_1)},\Sigma^2_{N_n(t_2)} - \Sigma^2_{N_n(t_1)}\}. \]
		Now since $\Sigma^2_{j,j+1}=\E_m[(\bar{v}_j \circ \cal{T}^j)^2]\le |\bar v_j|_\infty^2$, Lemma~\ref{lem:covar_approx_block_sum} implies that $|\Sigma^2_{j+1}-\Sigma^2_j|\lesssim (j+1)^{2-\gamma}+\log (j+1)$ for all $j\ge 0$. Since $N_n(t_2)\le n$, it follows that \[|\Sigma^2_{N_n(t_2)}-\Sigma^2_{N_n(t_1)} - (t_2 - t_1)\Sigma_n^2|\lesssim n^{2-\gamma}+\log n.\]
		Similarly, another application of Lemma~\ref{lem:covar_approx_block_sum} yields that 
		\[ |\Sigma^2_{N_n(t_2)}-\Sigma^2_{N_n(t_1)} - \Sigma^2_{N_n(t_1),N_n(t_2)}|\lesssim n^{2-\gamma}+\log n \]
		and so $|(t_2 - t_1)\Sigma_n^2 - \Sigma^2_{N_n(t_1),N_n(t_2)}|\lesssim n^{2-\gamma}+\log n$. The result follows since $\E_m[(W_n(t_2)-W_n(t_1))^2]=\Sigma_n^{-2}\Sigma^2_{N_n(t_1),N_n(t_2)}$.
	\end{proof}
    
\begin{prop}\label{prop:X_i-var-approx}
There is a constant $C>0$ such that for all $n\ge 1$, we have
\[|\E_m[X_i^2]-1/k|\le C(\rho^2+\eta+\Sigma_n^{-2}(n^{2-\gamma}+\log n)).\]
\end{prop}
\begin{proof}
Recall that $k=[1/(\rho+\eta)]$, so 
\[0\le \frac{1}{\rho}-k\le \frac{1}{\rho}-\frac{1}{\rho+\eta}+1\le \frac{\eta}{\rho^2}+1.\] 
Since $1/k\le 4\rho$ by~\eqref{eq:k-asymp}, it follows that
\[ \bigg|\frac1k -(s_i-t_{i-1})\bigg|=\bigg|\frac1k -\rho\bigg|=\frac{\rho}{k}\Big|\frac1\rho-k\Big|\le 4(\rho^2+ \eta). \]
On the other hand, $|\E_m[X_i^2]-(s_i-t_{i-1})|\le C\Sigma_n^{-2}(n^{2-\gamma}+\log n)$ by Proposition~\ref{prop:block_var_approx}. The result follows by combining these estimates.
\end{proof}

\begin{lem}\label{lem:small_block_max_moment}
Let $\eps>0$. There exists a constant $C>0$ such that for all $n\ge 1$, we have
\[ \frac{1}{\Sigma_n}\bigg\|\max_{m\le n}\Big|\sum_{i=0:i\in \cal J}^{m-1}\bar v_i \circ \cal T^i\Big| \bigg\|_{\gamma+1}\le Cn^\eps\Big[k^{\frac12}\eta^{\frac12}+\rho^{\frac12} +\Sigma_n^{-1}\big(k^{\frac12}(n^{1-\frac{\gamma}{2}}+\log n)+ n^{\frac{1}{\gamma+1}}\big)\Big]. \]
\end{lem}

\begin{proof}
Write $\cal J = \{\ell_0, \dots,\ell_{d-1} \}$ where $\ell_0<\ell_1<\dots<\ell_{d-1}$ so $d\le n$. For $i\ge d$ define $\ell_i =\ell_{d-1}+ i$. Note that 
\begin{align}\label{eq:Y_i_cum_sum}
    \max_{m\le n}\Big|\sum_{i=0:i\in \cal J}^{m-1}\bar v_i \circ \cal T^i\Big|= \max_{u\le d}\Big|\sum_{i=0}^{u-1}\bar v_{\ell_i}\circ \cal T^{\ell_i}\Big|.
\end{align}
As noted in Remark~\ref{rem:fcb}(ii), the sequence $\{\cal T^{\ell_i}\}_{i\ge 0}$ satisfies the Functional Correlation Bound with rate $n^{-\gamma}$ and the same constant $C_{\fcb}$ as the sequence $\{\cal T^n\}_{n\ge 0}$. Hence by Proposition~\ref{prop:moment_bds_fcb},
\begin{align}\label{eq:small_block_moment_bd}
    \bigg\|\max_{u\le d}\Big|\sum_{i=0}^{u-1}\bar v_{\ell_i}\circ \cal T^{\ell_i}\Big| \bigg\|_{\gamma+1}^{\gamma+1} \lesssim n^\eps \bigg(\Big(\Var \sum_{i=0}^{d-1}\bar v_{\ell_i}\circ \cal T^{\ell_i}\Big)^{\frac{\gamma+1}{2}}+n\bigg).
\end{align}
Recall that $\cal J = \bigcup_{i=1}^{k}\cal J_i$ so by applying Lemma~\ref{lem:covar_approx_block_sum} (again with $\{\cal T^{\ell_{i+1}}\}$ in place of $\{\cal T^n\}$ and $\{v_{\ell_i}\}$ in place of $\{v_n\}$), we obtain that 
\begin{align}\label{eq:var_small_blocks}
   \nonumber \Var \bigg(\sum_{i=0}^{d-1}\bar v_{\ell_i}\circ \cal T^{\ell_i}\bigg)&\lesssim \sum_{i=1}^k \Var\Big(\sum_{j\in \cal J_i}\bar v_j \circ \cal T^j\Big)+C\sum_{i=1}^{k}\big(|\cal J_i|^{2-\gamma}+\log|\cal J_i|\big)\\ 
    &\lesssim \sum_{i=1}^k \Sigma_n^2\Var Y_i + k^{\gamma-1}n^{2-\gamma}+k\log n.
\end{align}
Indeed, if $1<\gamma<2$, then by H\"older's inequality,
\[
    \sum_{i=1}^{k}|\cal J_i|^{2-\gamma}
    \le k^{\gamma-1}\Big(\sum_{i=1}^k |\cal J_i|\Big)^{2-\gamma}
    \le k^{\gamma-1}n^{2-\gamma}.
\]
If $\gamma\ge2$, then $\sum_{i=1}^{k}|\cal J_i|^{2-\gamma}\le k\le k\log n$ since $\eta>1/n$ and $\rho\le1/4$ imply $n\ge5$. Moreover, $\sum_{i=1}^{k}\log|\cal J_i|\le k\log n$.
Recall that $Y_i = W_n(t_i)-W_n(s_i)$ by~\eqref{eq:blocks_W_n}. Moreover, $t_i - s_i = \eta$ for $1\le i\le k-1$ and $t_k-s_k\le 2(\rho+\eta)$. Hence by Proposition~\ref{prop:block_var_approx} we obtain that 
\begin{align*}
    \sum_{i=1}^k \Var Y_i &\lesssim (k-1)\eta+2(\rho+\eta)+k\Sigma_n^{-2}(n^{2-\gamma}+\log n)\\
    &\lesssim k\eta+\rho+k\Sigma_n^{-2}(n^{2-\gamma}+\log n).
\end{align*}
Finally, by combining this bound with~\eqref{eq:Y_i_cum_sum}--\eqref{eq:var_small_blocks} we conclude that
\begin{align*}
    &\frac{1}{\Sigma_n}\bigg\|\max_{m\le n}\bigg|\sum_{i=0:i\in \cal J}^{m-1}\bar v_i \circ \cal T^i\Big| \bigg\|_{\gamma+1} \\
    \lesssim &\Sigma_n^{-1} n^\eps \bigg[\Big(\Sigma_n^2\big(k\eta+\rho+k\Sigma_n^{-2}(n^{2-\gamma}+\log n)\big)+k^{\gamma-1}n^{2-\gamma}+k\log n\Big)^{\frac{\gamma+1}{2}}+n\bigg]^{\frac{1}{\gamma+1}}\\
    \lesssim &
    n^\eps\Big[k^{\frac12}\eta^{\frac12}+\rho^{\frac12} +\Sigma_n^{-1}\big(k^{\frac12}(n^{1-\frac{\gamma}{2}}+\log n)+ n^{\frac{1}{\gamma+1}}\big)\Big],
\end{align*}
as required.
\end{proof}

\begin{lem}\label{lem:big_block_max_moment}
Let $\eps>0$. There exists a constant $C>0$ such that for all $n\ge 1$, we have
\[
\Big\|\max_{1\le i\le k}\big|U_i\big|\Big\|_{\gamma+1}\le Cn^\eps \big[k^{\frac{1}{\gamma+1}}\big(\rho^{\frac12}+\Sigma_n^{-1}(n^{1-\gamma/2}+\log n)\big)+\Sigma_n^{-1}n^{\frac{1}{\gamma+1}}\big].
\]
\end{lem}
\begin{proof}
We first observe that
\begin{align*}
    \Big\|\max_{1\le i\le k}\big|U_i\big|\Big\|_{\gamma+1}=\E_m\bigg[\max_{1\le i\le k}\big|U_i\big|^{\gamma+1}\bigg]^{\frac{1}{\gamma+1}}\le \Big(\sum_{i=1}^{k}\E_m\big[|U_i|^{\gamma+1}\big]\Big)^{\frac{1}{\gamma+1}}.
\end{align*}
Since the sum of the lengths of the big blocks is less than $n$, by Propositions~\ref{prop:moment_bds_fcb} and~\ref{prop:X_i-var-approx}, we have
\begin{align}
\Big(\sum_{i=1}^{k}\E_m\big[|U_i|^{\gamma+1}\big]\Big)^{\frac{1}{\gamma+1}}\lesssim & n^\eps\Big(\sum_{i=1}^{k}\big(\E_m[|X_i|^{2}]\big)^{\frac{\gamma+1}{2}}+\Sigma_n^{-(\gamma+1)}n\Big)^{\frac{1}{\gamma+1}} \nonumber\\
\lesssim &n^\eps\Big(k\big(\tfrac{1}{k}+\rho^2+\eta+\Sigma_n^{-2}(n^{2-\gamma}+\log n)\big)^{\frac{\gamma+1}{2}}+\Sigma_n^{-(\gamma+1)}n\Big)^{\frac{1}{\gamma+1}} \nonumber\\
\lesssim & n^\eps \Big(k^{\frac{1}{\gamma+1}}(\rho^{1/2}+\eta^{1/2}+\Sigma_n^{-1}(n^{1-\gamma/2}+\log n))+\Sigma_n^{-1}n^{\frac{1}{\gamma+1}}\Big), \label{eq:sum_pth_moments_big_blocks}
\end{align}
where we have used the fact that $1/k\le 4\rho$ in the last line. Since $\eta<\rho$, we obtain the result as required.
\end{proof}
\medskip
{\bf Step 2. Let $\widehat X_1, \widehat X_2,\ldots, \widehat X_k$ be independent random variables with $\widehat X_i=_d X_i$. Define $\widehat W_n(t)=\sum_{i=1}^{[kt]}\widehat X_i$.
Bound $\mathcal{W}_{1}(\widetilde W_n,\widehat W_n)$ by using the Functional Correlation Bound.}
\vskip 2mm

			\begin{lem}\label{xnw}
				There exists a constant $C>0$ such that $\mathcal{W}_1(\widetilde{W}_n,\widehat{W}_n)\le C k^2 n\Sigma_n^{-(1+2\gamma)}\eta^{-\gamma}$ for all $n\ge 1$.
			\end{lem}

\begin{proof}
We follow the argument in \cite[Lemma~4.3]{F25}. By the Kantorovich-Rubinstein theorem (cf.~\cite[Theorem 3.2.7]{bogachev_weak_conv}), it suffices to show that for all $n\ge 1$, we have
\[
\big|\E_m[G(\widetilde W_n)]-\E[G(\widehat W_n)]\big|\le C k^2 n\Sigma_n^{-(1+2\gamma)}\eta^{-\gamma}
\]
for any Lipschitz function $G:D[0,1]\to \R$ such that $\Lip(G)\le 1$ and $G(0)=0$.

Now, for $1\le i\le k$, we can write $X_i(x)=\Phi_i(\cal T^{a_i}x,\ldots,\cal T^{b_i}x)$, where $a_i=N_n(t_{i-1})$, $b_i=N_n(s_i)-1$ and 
\[
\Phi_i(y_{a_i},\ldots,y_{b_i})=\frac{1}{\Sigma_n}\sum_{j=a_i}^{b_i}\bar v_j(y_j).
\]

Let $R=\max_{1\le i\le k}|\Phi_i|_\infty$ and define $\pi_k:[-R,R]^k\to D[0,1]$ by
\[
\pi_k(x_1,\ldots,x_k)(t)=\sum_{i=1}^{[kt]}x_i, \text{~for~} t\in[0,1].
\]
Then $\widetilde W_n=\pi_k(X_1,\ldots,X_k)$, $\widehat W_n=\pi_k(\widehat X_1,\ldots,\widehat X_k)$. Observe that $G\circ \pi_k$ is Lipschitz continuous. By applying Proposition~\ref{lip}
\footnote{To apply Proposition~\ref{lip}, we omit any empty big blocks; this has no effect on the argument since the corresponding $X_i$ are identically zero.}with $\varphi=G\circ \pi_k$, we have
\[
\big|\E_m[G(\widetilde W_n)]-\E[G(\widehat W_n)]\big|\le A,
\]
where $A=C_{\fcb}\sum_{r=1}^{k-1}(a_{r+1}-b_r)^{-\gamma}\Big(|G\circ \pi_k|_\infty+\Lip(G\circ \pi_k)\sum_{i=1}^{k}\sum_{j=1}^{b_i - a_i + 1}\Lip(\Phi_i;j)\Big)$.

Since $\Sigma_{N_n(s_i)}^2\le s_i\Sigma_n^2$ and $\Sigma_{N_n(t_i)+1}^2> t_i \Sigma_n^2$, by Lemma~\ref{lem:var_diff}, we have
\begin{align*}
    a_{i+1} - b_i = N_n(t_i) - N_n(s_i)+1\gtrsim \Sigma^2_{N_n(t_i)+1}-\Sigma^2_{N_n(s_i)}> (t_i - s_i)\Sigma_n^2=\eta \Sigma_n^2.
\end{align*}
It follows that
\[
\sum_{r=1}^{k-1}(a_{r+1}-b_r)^{-\gamma}\lesssim k(\Sigma_n^2 \eta)^{-\gamma}.
\]
Also, since $b_i - a_i + 1\le n$ for all $i$, we have $R\lesssim \Sigma_n^{-1}n$. Since $G(0)=0$ by assumption, it follows that $|G\circ \pi_k|_\infty\le \Lip(G)|\pi_k|_\infty\le kR \,\Lip(G)\lesssim \frac{kn}{\Sigma_n}$. Moreover, $\Lip(\pi_k)\le 1$ and thus
\[
\Lip(G\circ \pi_k)\sum_{i=1}^{k}\sum_{j=1}^{b_i - a_i +1}\Lip(\Phi_i;j)\le n\max_{i,j} \Lip(\Phi_i;j)\le \tfrac{n}{\Sigma_n}\sup_i \Lip(\bar v_i).
\]
Combining the above estimates, we obtain the result.
\end{proof}

{\bf Step 3. Estimation of the Wasserstein convergence rate between $\widehat W_n$ and $B$, where $B$ is a standard Brownian motion.}

\vskip3mm

We use the following strong approximation result for independent random vectors due to Sakhanenko~\cite{S85}.
\begin{prop}[{\cite{S85}}]\label{indep}
Let $\widehat X_1, \widehat X_2,\ldots, \widehat X_k$ be independent random variables with $\E[\widehat X_i]=0$, $i=1,\ldots, k$. Let $p>2$ and let $L_p=\sum_{i=1}^{k}\E[|\widehat X_i|^p]<\infty$. Then one can construct on a probability space a sequence of independent random variables $\widetilde X_1,\ldots,\widetilde X_k$ and a corresponding sequence of
independent Gaussian random variables $Z_1,\ldots, Z_k$ such that $\widetilde X_i=_d \widehat X_i$, $\E [Z_i]=0$, $\Var Z_i=\Var \widetilde X_i$, $i=1,\ldots, k$ and
\[
\E \Big[\max_{1\le m\le k}\big|\sum_{i=1}^{m}\widetilde X_i-\sum_{i=1}^{m}Z_i\big|^p\Big]\le Cp^{2p}L_p.
\]
\end{prop}

\begin{lem}\label{lem:Brownian_indep}
Let $\eps>0$. There is a constant $C>0$ such that for all $n\ge 1$,
\[
\cal W_{\gamma+1}(\widehat W_n, B)\le Cn^\eps\Big[
        k^{\frac{1}{\gamma+1}-\frac12}+k^{\frac12}\big(\eta^{\frac12}+\Sigma_n^{-1}(n^{1-\gamma/2}+\log n)\big)+\Sigma_n^{-1}n^{\frac{1}{\gamma+1}}
    \Big].
\]
\end{lem}

\begin{proof}
Let $p=\gamma+1$ and define $Y_n(t)=\sum_{i=1}^{[kt]}Z_i$. Then
\[
{\mathcal W}_{p}(\widehat W_n, B)\le {\mathcal W}_{p}(\widehat W_n, Y_n)+{\mathcal W}_{p}(Y_n,B)=:I_1+I_2.
\]
By Proposition~\ref{indep}, we can estimate that
\begin{align*}
&I_1:={\mathcal W}_{p}(\widehat W_n, Y_n)\\
\le& \Big\|\max_{1\le m\le k}\big|\sum_{i=1}^{m}\widetilde X_i-\sum_{i=1}^{m}Z_i\big|\Big\|_{p}
\le C\Big(\sum_{i=1}^{k}\E\big[|\widehat X_i|^{p}\big]\Big)^{1/p}.
\end{align*}
Now, by~\eqref{eq:sum_pth_moments_big_blocks} and $\rho\le \frac{1}{k}$,
\begin{align*}
    \Big(\sum_{i=1}^{k}\E\big[|\widehat X_i|^{p}\big]\Big)^{\frac1p}=\Big(\sum_{i=1}^{k}\E_m\big[|X_i|^{p}\big]\Big)^{\frac1p}\lesssim n^\eps \Big(k^{\frac1p}\big(k^{-\frac12}+\eta^{\frac12}+\Sigma_n^{-1}(n^{1-\gamma/2}+\log n)\big)+\Sigma_n^{-1}n^{\frac1p}\Big).
\end{align*}

Next, we aim to estimate $I_2$.  Write $\tau_j = \sum_{i=1}^j \E_m[X_i^2]=\sum_{i=1}^j \E[Z_i^2]$. Note that $(B(\tau_j)-B(\tau_{j-1}))_{j=1}^k$ are independent centered Gaussian random variables and $\E[(B(\tau_j)-B(\tau_{j-1}))^2]=\tau_j - \tau_{j-1}=\E[Z_j^2]$, so $Y_n =_d \widetilde{Y}_n$ where
\[ \widetilde{Y}_n(t)=B\big(\tau_{[kt]}\big)=\sum_{j=1}^{[kt]}(B(\tau_j)-B(\tau_{j-1})). 
\]
Let $q\ge p$. It follows that
\begin{align}
 \cal W_p(Y_n,B)&=\cal W_p(\widetilde Y_n, B)\le \bigg\|\sup_{t\in[0,1]}|\widetilde Y_n(t)-B(t)| \bigg\|_q\le \bigg\|\max_{0\le j\le k}\,\sup_{t\in[j/k,(j+1)/k)}|B(\tau_j)-B(t)| \bigg\|_q\nonumber\\
&\le (k+1)^{1/q}\max_{0\le j\le k}\bigg\|\sup_{t\in[j/k,(j+1)/k)}|B(\tau_j)-B(t)| \bigg\|_q\label{eq:brownian_initial_bd}.
\end{align}
Now by Proposition~\ref{prop:X_i-var-approx}, for $0\le j\le k$ we have
\begin{align*}
&\big\|B(\tau_j)-B(j/k)\big\|_q\lesssim |\tau_j-j/k|^{1/2}\\
=&\Big|\sum_{i=1}^j (\E_m[X_i^2]-1/k)\Big|^{1/2}\\
\lesssim & k^{1/2}\big(\rho+\eta^{1/2}+\Sigma_n^{-1}(n^{1-\gamma/2}+\log n)\big).
\end{align*}
Hence by the reflection principle we obtain that
\begin{align*}
\bigg\|\sup_{t\in[j/k,(j+1)/k)}|B(\tau_j)-B(t)| \bigg\|_q
&\le \big\|B(\tau_j)-B(j/k)\big\|_q+\bigg\|\sup_{t\in[j/k,(j+1)/k)}|B(j/k)-B(t)| \bigg\|_q\\
&\lesssim k^{1/2}\big(\rho+\eta^{1/2}+\Sigma_n^{-1}(n^{1-\gamma/2}+\log n)\big)+k^{-1/2}\\
&\lesssim k^{1/2}\big(\eta^{1/2}+\Sigma_n^{-1}(n^{1-\gamma/2}+\log n)\big)+k^{-1/2},
\end{align*}
where we have used that $k\le \frac{1}{\rho}$ so $k^{1/2}\rho\le k^{-1/2}$. Combining this bound with~\eqref{eq:brownian_initial_bd}, we obtain that
\[
I_2=\cal W_p(Y_n,B)\lesssim k^{\frac{1}{2}+\frac{1}{q}}\big(\eta^{\frac12}+\Sigma_n^{-1}(n^{1-\gamma/2}+\log n)\big)+k^{-\frac{1}{2}+\frac{1}{q}}.
\]
Let $q\ge 1/\eps$. Since $\frac{1}{4\rho}\le k\le \frac{1}{\rho}\le n$ by~\eqref{eq:k-asymp} we thus obtain that
\begin{align*}
    &I_1 + I_2 \\
    \lesssim &n^\eps \Big(k^{\frac1p}\big(k^{-\frac12}+\eta^{\frac12}+\Sigma_n^{-1}(n^{1-\gamma/2}+\log n)\big)+\Sigma_n^{-1}n^{\frac1p}\Big)\\
    &\qquad +k^{\frac{1}{2}+\eps}\big(\eta^{\frac12}+\Sigma_n^{-1}(n^{1-\gamma/2}+\log n)\big)+k^{-\frac{1}{2}+\eps}\\
    \lesssim& n^\eps\Big(k^{\frac1p-\frac12}+k^{\frac12}\big(\eta^{\frac12}+\Sigma_n^{-1}(n^{1-\gamma/2}+\log n)\big)+\Sigma_n^{-1}n^{\frac1p} \Big),
\end{align*}
as required.
\end{proof}

One last elementary lemma is needed in the proof of our main theorem.
	\begin{lem}\label{lem:opt_problem}
		Let $a,b,c>0$ and $K_1,K_2>0$. Set $M=\max\{K_1^{1/(a+b)},K_2^{1/(a+c)}\}$. Then there exists a constant $C>1$ depending only on $a,b,c$ such that 
		\begin{equation}\label{eq:opt_problem_bds}
			C^{-1}M^a\le \inf_{\rho>0}[\rho^a + K_1 \rho^{-b}+K_2\rho^{-c}]\le CM^a.
		\end{equation}
		Moreover, the infimum is attained at a point $\rho\in [C^{-1}M,CM]$.
	\end{lem}
    We do not use the lower bound in~\eqref{eq:opt_problem_bds}; we include it only to show that the upper bound is sharp.
	\begin{proof}
		For $\rho>0$ define $g(\rho)=\rho^a + K_1 \rho^{-b}+K_2\rho^{-c}$.
		 The unique stationary point of $g$ satisfies $a\rho_*^a = K_1 b \rho_*^{-b}+K_2 c\rho_*^{-c}$. Since $\max\{K_1 b\rho_*^{-b},K_2 c \rho_*^{-c}\}\le a\rho_*^a \le 2\max\{K_1 b\rho_*^{-b},K_2 c\rho_*^{-c}\}$, we have
		\[ \max\big\{\big(\tfrac{b}{a}K_1\big)^{\frac{1}{a+b}},\ \big(\tfrac{c}{a}K_2\big)^{\frac{1}{a+c}} \big\}\le \rho_*
		\le
		\max\big\{\big(\tfrac{2b}{a}K_1\big)^{\frac{1}{a+b}},\ \big(\tfrac{2c}{a}K_2\big)^{\frac{1}{a+c}}\big\}. \]
		That is, $C^{-1}M\le \rho_* \le CM$ where $C$ depends only on $a,b,c$. Since $C^{-a}M^a\le \rho_*^a \le C^a M^a$, the lower bound in~\eqref{eq:opt_problem_bds} follows after changing $C$ if necessary. Finally, the upper bound in~\eqref{eq:opt_problem_bds} follows by noting that \[K_1 \rho_*^{-b}\le K_1 C^{b}M^{-b}\le C^{b}K_1^{1-b/(a+b)}=C^{b}K_1^{a/(a+b)}\]
		and similarly $K_2 \rho_*^{-c}\le C^{c}K_2^{a/(a+c)}$.
    \end{proof}

We are now ready to prove our main theorem.

\begin{proof}[Proof of Theorem~\ref{thnon}]
We first deal with the trivial case where $\Sigma_n \le n^{\frac{1}{\gamma+1}}$. By Proposition~\ref{prop:moment_bds_fcb}, 
        \begin{equation}\label{eq:W_n-norm}
            \cal{W}_{\gamma+1}(0,W_n)\le\Big\|\Sigma_n^{-1}\max_{1\le j\le n}|S_j|\Big\|_{\gamma+1}\lesssim \Sigma_n^{-1}n^\eps (\Sigma_n+n^{\frac{1}{\gamma+1}}).
        \end{equation}
Since $\cal W_{\gamma+1}(0,B)\lesssim 1\le \Sigma_n^{-1}n^{\frac{1}{\gamma+1}}$, it follows that $\cal W_{\gamma+1}(W_n,B)\lesssim n^\eps \Sigma_n^{-1}n^{\frac{1}{\gamma+1}}$, which concludes the proof of the theorem in this case.

From now on we assume $\Sigma_n \ge n^{\frac{1}{\gamma+1}}$. Let $1\le r<\gamma+1$. By the triangle inequality,
\begin{equation}\label{eq:main_triangle}
    \mathcal W_r(W_n,B)
    \le \mathcal W_r(W_n,\widetilde W_n)
       +\mathcal W_r(\widetilde W_n,\widehat W_n)
       +\mathcal W_r(\widehat W_n,B).
\end{equation}
		We estimate these three terms separately. First, by Lemmas~\ref{lem:small_block_max_moment} and~\ref{lem:big_block_max_moment}, for every $\eps>0$,
		\begin{align*}%
			\mathcal W_{\gamma+1}(W_n,\widetilde W_n)
			&\lesssim n^\eps\bigg[
			k^{\frac12}\eta^{\frac12}+\rho^{\frac12} +\Sigma_n^{-1}\big(k^{\frac12}(n^{1-\frac{\gamma}{2}}+\log n)+ n^{\frac{1}{\gamma+1}}\big) \notag \\
			&\qquad\qquad\qquad
			+k^{\frac{1}{\gamma+1}}\big(\rho^{\frac12}+\Sigma_n^{-1}(n^{1-\frac{\gamma}{2}}+\log n)\big)+\Sigma_n^{-1}n^{\frac{1}{\gamma+1}}\bigg] \notag\\
			&\lesssim 
            n^\eps\bigg[k^{\frac{1}{\gamma+1}}\rho^{\frac12}+ k^{\frac12}\big(\eta^{\frac12}+\Sigma_n^{-1}(n^{1-\frac{\gamma}{2}}+\log n)\big)
            + \Sigma_n^{-1}n^{\frac{1}{\gamma+1}}  \bigg]. 
		\end{align*}
		By Lemma~\ref{lem:Brownian_indep}, we have
		\begin{align*}%
			\mathcal W_{\gamma+1}(\widehat W_n,B) 
			\lesssim n^\eps\bigg[k^{\frac{1}{\gamma+1}-\frac12}+k^{\frac12}\big(\eta^{\frac12}+\Sigma_n^{-1}(n^{1-\frac{\gamma}{2}}+\log n)\big) + \Sigma_n^{-1} n^{\frac{1}{\gamma+1}} 
			\bigg].
		\end{align*}
        Recall that $\frac{1}{4\rho}\le k\le \frac{1}{\rho}$ by~\eqref{eq:k-asymp}, so it follows that $\cal W_{\gamma+1}(W_n,\widetilde W_n) + \cal W_{\gamma+1}(\widehat W_n,B)\lesssim n^\eps \Delta_1$, where 
        \begin{equation*}
            \Delta_1 = \rho^{\frac12-\frac{1}{\gamma+1}}
				+(\tfrac{\eta}{\rho})^{\frac12}+\Sigma_n^{-1}(n^{1-\frac{\gamma}{2}}+\log n)\rho^{-\frac12}+\Sigma_n^{-1}n^{\frac{1}{\gamma+1}}.
        \end{equation*}
		
		By Lemma~\ref{xnw}, $\mathcal W_1(\widetilde W_n,\widehat W_n)\lesssim D_n$, where
		\[
		D_n=k^2 n\Sigma_n^{-(1+2\gamma)}\eta^{-\gamma}.
		\]
		We now use this to obtain a bound for $\cal W_r(\widetilde W_n,\widehat W_n)$ by interpolation. Indeed, by Lemma~3.1 in \cite{F25}, we know that if
		$1\le r<p$ and $U,V$ are random elements of $D[0,1]$ with finite $p$-th moments,
		then
		\begin{equation}\label{eq:wasserstein_interpolation_inproof}
			\mathcal W_r(U,V)
			\le \mathcal W_1(U,V)^{\frac{p-r}{r(p-1)}}
			\Big(\mathcal W_p(U,0)+\mathcal W_p(0,V)\Big)^{1 - \frac{p-r}{r(p-1)}}.
		\end{equation}
		
		Since $\Sigma_n\ge n^{\frac{1}{\gamma+1}}$,~\eqref{eq:W_n-norm} implies that $\cal W_{\gamma+1}(0,W_n)=O(n^\eps)$.
		It follows that
			\[\cal W_{\gamma+1}(0,\widetilde W_n)+\cal W_{\gamma+1}(0,\widehat W_n)\lesssim n^\eps+ \cal W_{\gamma+1}(W_n,\widetilde W_n) + \cal W_{\gamma+1}(\widehat W_n,B)\lesssim n^\eps(1+\Delta_1).\]
			Consequently, by applying~\eqref{eq:wasserstein_interpolation_inproof} with $p=\gamma+1$, it follows that
			\begin{equation}\label{eq:main_step3_bound}
				\mathcal W_r(\widetilde W_n,\widehat W_n)
				\lesssim n^\eps D_n^{\kappa}(1+\Delta_1)^{1-\kappa}=:n^\eps \Delta_2 (1+\Delta_1)^{1-\kappa},
			\end{equation}
			where $\kappa=\frac{\gamma+1-r}{r\gamma}$. Hence we obtain $\cal W_r(W_n,B)\lesssim n^\eps (\Delta_1 + \Delta_2(1+\Delta_1)^{1-\kappa}
            )$. We claim that we can choose $\eta$, $\rho$ so that
			\begin{equation}\label{eq:Delta_1_2_claimed_bd}
				\Delta_1 + \Delta_2 \lesssim  (n\Sigma_n^{-(1+2\gamma)})^{\ell(\gamma,r)} + \big[\Sigma_n^{-1}(n^{1-\frac{\gamma}{2}}+\log n)\big]^{\frac{\gamma-1}{2\gamma}} + n^{\frac{1}{\gamma+1}}\Sigma_n^{-1}.
			\end{equation}
			Since $\Sigma_n\ge n^{\frac{1}{\gamma+1}}$, it follows that $\Delta_1 + \Delta_2 = O(1)$ and so $\cal W_r(W_n,B)\lesssim n^\eps(\Delta_1 + \Delta_2)$. The bounds stated in the theorem follow by applying~\eqref{eq:Delta_1_2_claimed_bd} again. Indeed, for $1<\gamma<2$ the stated bound follows immediately since $\log n \lesssim n^{1-\frac{\gamma}{2}}$.

            Let $\gamma\ge 2$. Since $1\le r<\gamma+1$, a straightforward calculation shows that $(1+2\gamma)\ell(\gamma,r)\le (\gamma-1)/(2\gamma)$ and it follows that \[\big[\Sigma_n^{-1}(n^{1-\frac{\gamma}{2}}+\log n)\big]^{\frac{\gamma-1}{2\gamma}}\lesssim n^{\ell(\gamma,r)}\Sigma_n^{-\frac{\gamma-1}{2\gamma}}\le (n\Sigma_n^{-(1+2\gamma)})^{\ell(\gamma,r)},\] which concludes the proof in this case.

		It remains to prove the claim~\eqref{eq:Delta_1_2_claimed_bd}. Note that
		\begin{equation}
			\Delta_1 + \Delta_2 \lesssim T_1 + T_2 + T_3 + T_4 + n^{\frac{1}{\gamma+1}}\Sigma_n^{-1},
		\end{equation}	
		where 
		\begin{align*} 
			T_1 = \Big(\rho^{-2} n\Sigma_n^{-(1+2\gamma)}\eta^{-\gamma}\Big)^{\kappa}, \quad
			 T_2 = (\tfrac{\eta}{\rho})^{1/2}, \quad 
			 T_3 = \rho^{\frac12-\frac{1}{\gamma+1}}, \quad
			 T_4 = \Sigma_n^{-1}(n^{1-\frac{\gamma}{2}}+\log n)\rho^{-\frac12}.
		\end{align*}
		We now choose $\eta,\rho$ so as to minimize $T_1 + T_2 + T_3 + T_4$. Note that $T_1$ is decreasing in $\eta$ and $T_2$ is increasing in $\eta$, whilst $T_3$ and $T_4$ are independent of $\eta$. Hence we choose $\eta$ by setting $T_1 = T_2$, which yields \[\eta=(n\Sigma_n^{-(1+2\gamma)})^{\frac{2\kappa}{1+2\kappa \gamma}}\rho^{\frac{1-4\kappa}{1+2\kappa \gamma}}.\] 
		We thus obtain 
		\begin{align*}
			T_1 + T_2 + T_3 + T_4 = \rho^a+
			K_1 \rho^{-b}+
			K_2 \rho^{-\frac12},
		\end{align*}
		where 
		\begin{align*}
			a&=\frac{1}{2}-\frac{1}{\gamma+1},\qquad b=\frac{1}{2}\bigg(1-\frac{1-4\kappa}{1+2\kappa\gamma}\bigg),\\
			K_1 &= 2(n\Sigma_n^{-(1+2\gamma)})^{\frac{\kappa}{1+2\kappa \gamma}}, \qquad
			K_2 = \Sigma_n^{-1}(n^{1-\frac{\gamma}{2}}+\log n).
		\end{align*}
		Hence by Lemma~\ref{lem:opt_problem}, we can choose $\rho$ such that
		\begin{equation}\label{eq:rho_asymp}
			\rho \asymp \max\{K_1^{\frac{1}{a+b}},K_2^{\frac{1}{a+1/2}}\}
		\end{equation}
		and
		\begin{align}
			T_1 + T_2 + T_3 + T_4 &\lesssim \max\{K_1^{\frac{a}{a+b}},K_2^{\frac{a}{a+1/2}}\}
			\label{eq:T_i_bds_prelim} \\
			&\lesssim \max\big\{(n\Sigma_n^{-(1+2\gamma)})^{\frac{\kappa}{1+2\kappa \gamma}\frac{a}{a+b}},\big[\Sigma_n^{-1}(n^{1-\frac{\gamma}{2}}+\log n)\big]^{\frac{a}{a+1/2}}\big\}.
			\notag
		\end{align}
		We have 
		\begin{align*}
			\frac{\kappa}{1+2\kappa\gamma}\cdot\frac{a}{a+b}&=\frac{\kappa(\gamma-1)}{4(\gamma+1)\gamma\kappa + \gamma + 4\kappa - 1} =  
			\ell(\gamma,r),\\
			\frac{a}{a+1/2}&=\frac{\gamma-1}{2\gamma}.
		\end{align*}
		This gives the bound stated in~\eqref{eq:Delta_1_2_claimed_bd}; finally, we need to check that the chosen parameters satisfy $\frac{1}{n}<\eta<\rho\le \frac{1}{4}$. Since $\Sigma_n\ge n^{\frac{1}{\gamma+1}}$, we have $K_1=o(1)$ and $K_2=o(1)$ and thus $\rho=o(1)$ by~\eqref{eq:rho_asymp}. Moreover,~\eqref{eq:T_i_bds_prelim} implies that $T_2=(\eta/\rho)^{1/2}=o(1)$ so $\eta=o(\rho)$. Applying~\eqref{eq:T_i_bds_prelim} again, $T_1^{1/\kappa}=\rho^{-2} n\Sigma_n^{-(1+2\gamma)}\eta^{-\gamma}=o(1)$. Since $\Sigma_n^2 \lesssim n$, it follows that
		\[ \Sigma_n^{-2\gamma}\eta^{-\gamma}\lesssim \rho^2 n^{-1} \Sigma_n\lesssim \rho^2 = o(1) \]
		and so $\Sigma_n^{-2}=o(\eta)$. By again using that $\Sigma_n^2 \lesssim n$, it follows that $\frac{1}{n}=o(\eta)$. Thus for $n$ large enough, $\frac{1}{n}<\eta<\rho\le \frac{1}{4}$, as required.
    \end{proof}

\appendix\label{sec:variance_growth_sharp}

\section{Sharpness of the variance growth assumption}

In this section, we provide a counterexample to show that the variance growth assumption in Theorem~\ref{thnon} is essentially sharp. 
Let $M=[0,1]$ and let $T:M\to M$ be the LSV map~\eqref{P-M} with parameter $0<\alpha<\frac12$. 
Let $m$ be the unique absolutely continuous invariant probability measure. 
By~\cite{Y99}, $T$ is modeled by a Young tower with $O(n^{-1/\alpha})$ tails. 
So by~\cite{F22}, for any $\eta\in(0,1]$, the Functional Correlation Bound holds with rate $n^{-\gamma}$ for separately $\eta$-H\"older functions, where $\gamma=1/\alpha-1$. 
As explained in Remark~\ref{rem:fcb}~(iii), by a change of metric, Theorem~\ref{thnon} also applies to sequences of observables $v_n:M\to \R$ such that $\sup_n \|v_n\|_\eta<\infty$, where $\|v\|_\eta=|v|_\infty+|v|_\eta$ and $|v|_\eta=\sup_{x\ne y}|v(x)-v(y)|/|x-y|^\eta$.

Let $0<\eta<\alpha$, $r>1$ and set $c_j =(j+1)^{-r}$. Let $\psi:[0,\infty)\to [0,1]$ be $\eta$-H\"older and supported on $[0,1]$ with $\psi|_{[0,\frac12]}=1$. Define $v_j = c_j^\eta \psi(x/c_j)$. Then $|v_n|_\eta=|\psi|_\eta$ and $|v_n|_\infty=c_n^\eta$. In particular, $\sup_n \|v_n\|_\eta<\infty$.
\begin{lem}
    Let $1<r<(1-\alpha+\eta)^{-1}$ and $\kappa=\frac12(1-r(1-2\alpha+2\eta))$. Then $\Sigma_n\gtrsim n^{\kappa}$ and $S_n/\Sigma_n \to 0$ in $L^1$.
\end{lem}
We see that in this example the CLT fails and hence so does the WIP. 
Note that as $\eta\downarrow 0$ and $r\downarrow 1 $, we have $\kappa\uparrow \alpha=1/(\gamma+1)$. 
Thus the variance growth assumption that $\lim_{n\to\infty}\Sigma_n/n^\kappa=\infty$ for some $\kappa>1/(\gamma+1)$ in Theorem~\ref{thnon} is essentially sharp.
\begin{proof}
    By~\cite{thaler_80}, $\rmd m/\rmd \Leb\asymp x^{-\alpha}$ so $m([0,x])\asymp x^{1-\alpha}$.
    Since $v_j$ is supported on $[0,c_j]$, we have $\E_m[v_j]\le c_j^\eta m([0,c_j])\asymp c_j^{1-\alpha+\eta}$. Let $Z_n = \sum_{i=0}^{n-1}v_i\circ T^i$ so that $S_n=Z_n -\E_m[Z_n]$. Since $Z_n\ge 0$ and $-r(1-\alpha+\eta)\ne -1$, it follows that 
    \begin{align}\label{eq:L1_upper_bd}
        \|S_n\|_1 \le 2\|Z_n\|_1 =2\E_m[Z_n]\lesssim \sum_{j=0}^{n-1}c_j^{1-\alpha+\eta}=\sum_{j=0}^{n-1}(j+1)^{-r(1-\alpha+\eta)}\lesssim n^{\lambda},
    \end{align}
    where $\lambda=1-r(1-\alpha+\eta)$.

    Before proceeding further, we note that since $\alpha<\frac12$, we have $\alpha r<(1-\alpha+\eta)r<1$.
    Let $a_k = T|_{[0,1/2)}^{-k}(\tfrac12)$. Then by~\cite[Section 6.2]{Y99}, $a_k\asymp k^{-1/\alpha}$ so in particular there exists $L>0$ such that $y_{[Li^{\alpha r}]}\le \frac12 c_{2i}$ for all $i$. Let $i\le j\le i+i^{\alpha r}$. Then since $\alpha r<1$, we have $j\le 2i$ so
    \[ T^{j-i}[0,a_{[(L+2)i^{\alpha r}]}]=[0,a_{[(L+2)i^{\alpha r}]-(j-i)}]\subseteq [0,a_{[Li^{\alpha r}]}]\subseteq [0,\tfrac12 c_j]. \]
    Thus $v_i v_j\circ T^{j-i}=c_i^\eta c_j^\eta$ on $[0,a_{[(L+2)i^{\alpha r}]}]$ and it follows that for $i>0$,
    \[ \int v_i v_j\circ T^{j-i}\rmd m\gtrsim c_i^\eta c_j^\eta m([0,a_{[(L+2)i^{\alpha r}]}])\gtrsim c_i^\eta c_j^\eta \big({(i^{\alpha r})}^{-\frac{1}{\alpha}}\big)^{(1-\alpha)}\gtrsim i^{-r(2\eta+1-\alpha)}.\]
    For $i< n/2$,  we have $i+i^{\alpha r}\le 2i \le n$. Since $v_i\circ T^i v_j\circ T^j\ge 0$ for all $i,j$, it follows that
    \begin{align*}
        \E_m[Z_n^2]&=\sum_{i=0}^{n-1}\sum_{j=0}^{n-1}\int v_i\circ T^i v_j\circ T^j \rmd m\ge \sum_{0<i<n/2}\sum_{j=i}^{i+[i^{\alpha r}]}\int v_i \circ T^i v_j\circ T^j \rmd m\\
        &= \sum_{0<i<n/2}\sum_{j=i}^{i+[i^{\alpha r}]}\int v_i v_j\circ T^{j-i} \rmd m\\
        &\gtrsim \sum_{0<i<n/2}[i^{\alpha r}]i^{-r(2\eta+1-\alpha)}\gtrsim n^{1-r(1-2\alpha+2\eta)}=n^{2\kappa}.
    \end{align*}
    Finally, by~\eqref{eq:L1_upper_bd}, since $\lambda<\frac12+\frac12 r-r(1-\alpha+\eta)=\kappa$ we have $\Sigma_n^2=\E_m[Z_n^2]-\E_m[Z_n]^2\gtrsim n^{2\kappa}$ and moreover $\|S_n\|_1/\Sigma_n\lesssim n^{\lambda-\kappa}\to 0$.
\end{proof}

\section*{AI Usage Statement}
We used an AI tool (ChatGPT Sol 6.0) to help construct the counterexample presented in the appendix. Its sole purpose is to show that the variance growth assumption in the main theorem is essentially sharp. 
All other arguments given in the article were devised by the authors.
\section*{Acknowledgements}
NFV's research was supported by an Australian Research Council Laureate Fellowship (FL230100088).
ZW's research was supported by NSFC (Grant 12501239).

\bibliography{ref}
\bibliographystyle{alpha}

\end{document}